%% file: tLap09_12.tex
\documentclass[10pt]{article}
\usepackage{amscd,amssymb,amsmath,latexsym,enumerate,subfigure,amsthm,cite}
\usepackage[all]{xy}
\usepackage[mathscr]{euscript}
\usepackage{epsfig}
\usepackage{epsfig,psfrag}
\usepackage{float}
\usepackage{graphicx}

\usepackage{a4wide}

\usepackage{color,epsfig}

\title{Transverse Laplacians for Substitution Tilings\thanks{Work supported by
the NSF grants no.~DMS-0300398 and no.~DMS-0600956.}}

\author{Antoine Julien$^1$, Jean Savinien$^{2,3}$\\
 {\small $^1$ Institut Camille Jordan, Universit\'e Lyon I, France} \\
 {\small $^2$ Georgia Institute of Technology, Atlanta GA} \\
 {\small $^3$ SFB 701, Universit\"at Bielefeld, Germany.}
 }

 \date{ }

\theoremstyle{plain}
\newtheorem{theo}{Theorem}[section]
\newtheorem{proposi}[theo]{Proposition}
\newtheorem{lemma}[theo]{Lemma}
\newtheorem{coro}[theo]{Corollary}
\newtheorem{hypo}[theo]{Hypothesis}

\theoremstyle{definition}
\newtheorem{defini}[theo]{Definition}
\newtheorem{exam}[theo]{Example}

\theoremstyle{remark}
\newtheorem{rem}[theo]{Remark}
\newtheorem*{acks}{Acknowledgements}

\newtheoremstyle{citing}
  {\topsep}
  {\topsep}
  {\itshape}
  {}
  {\bfseries}
  {.}
  {.5em}
  {\thmnote{#3}}

\theoremstyle{citing}
\newtheorem*{varthm}{}

\newcommand{\Aa}{{\mathcal A}}
\newcommand{\Bb}{{\mathcal B}}

\newcommand{\Ee}{{\mathcal E}}

\newcommand{\Hh}{{\mathcal H}}

\newcommand{\Ll}{{\mathcal L}}

\newcommand{\Nn}{{\mathcal N}}
\newcommand{\Oo}{{\mathcal O}}
\newcommand{\Pp}{{\mathcal P}}

\newcommand{\Vv}{{\mathcal V}}

\newcommand{\CM}{{\mathbb C}}

\newcommand{\NM}{{\mathbb N}}
\newcommand{\QM}{{\mathbb Q}}

\newcommand{\RM}{{\mathbb R}}
\newcommand{\RMbar}{\overline{\mathbb R}}

\newcommand{\ZM}{{\mathbb Z}}

\newcommand{\TR}{{\rm Tr\,}}                       
\newcommand{\tra}{\mbox{\sc t}}                    
\newcommand{\Cs}{$C^{\ast}$-algebra }              
\newcommand{\Css}{$C^{\ast}$-algebras }            
\newcommand{\Sp}{\mbox{\rm Sp}}                    

\newcommand{\diam}{\mbox{\rm diam}}                

\newcommand{\pf}{\Lambda_{PF}}         
\newcommand{\ext}{\text{\rm ext}_1}    

\DeclareMathOperator{\Supp}{Supp}
\DeclareMathOperator{\Ext}{Ext}
\newcommand{\eps}{\varepsilon}
\DeclareMathOperator{\Card}{Card}

\newcommand{\Emod}{\mathcal{E}}
\newcommand{\Etot}{\mathcal{E}_\textrm{tot}}

\newcommand{\sca}[2]{\left\langle #1, #2\right\rangle} 

\begin{document}

\maketitle

\begin{abstract}
\noindent
Pearson and Bellissard recently built a {\em spectral triple} ---~the data of
Riemanian noncommutative geometry~--- for ultrametric Cantor sets.
They derived a family of Laplace--Beltrami like operators on those sets.
Motivated by the applications to specific examples, we revisit their work for
the transversals of tiling spaces, which are particular self-similar Cantor
sets.
We use {\em Bratteli diagrams} to encode the self-similarity, and {\em
Cuntz--Krieger algebras} to implement it.
We show that the abscissa of convergence of the $\zeta$-function of the spectral triple
gives indications on the exponent of complexity of the tiling.
We determine completely the spectrum of the Laplace--Beltrami operators, give
an explicit method of calculation for their eigenvalues, compute their Weyl
asymptotics, and a Seeley equivalent for their heat kernels.
\end{abstract}


\tableofcontents


\section{Introduction and summary of the results}
\label{tLap09.sect-intro}

In a recent article \cite{PB09}, Pearson and Bellissard defined a spectral
triple ---~the data of Riemanian noncommutative geometry (NCG) \cite{Co94}~---
for ultrametric Cantor sets.
They used a construction due to Michon \cite{Mich85}: any ultrametric Cantor set
$(C,d)$ can be represented isometrically as the set of infinite paths on a
weighted rooted tree.
The tree defines the topology, and the weights
encode the distance.
The spectral triple is then given in terms of combinatorial data on the tree.

With this spectral triple, they could define several objects,
including a $\zeta$-function, a measure $\mu$, and a one-parameter family of
operators on $L^2(C,\mu)$, which were interpreted as Laplace--Beltrami
operators.
They showed the abscissa of convergence of the $\zeta$-function to be a fractal
dimension of the Cantor set (the upper box dimension).

The goal of Pearson and Bellissard was to build a spectral triple for the
transversals of tiling dynamical systems. This opened a new, geometrical
approach to the theory of tilings.
Until now, all the operator-algebraic machinery used to study tilings and
tiling spaces was coming from noncommutative \emph{topology}.
Striking application of noncommutative topology were the study of the $K$-theory
for tiling \Css \cite{Bel82,BBG92,Kel95,kel97}, and namely
\cite{AP98,For97,DHS99} for computations applied to substitution tilings.
A follow-up of this study was the gap-labeling theorems for Schr\"oedinger
operators \cite{vEl94,BKL01,BBG06,BO03,KP03} 
---~a problem already appearing in some of the previously cited articles.
Other problems include cyclic cohomology and index theorems (and applications to
the quantum Hall effect)
\cite{BESB94,KRSB02,KSBa04,KSBb04,Kel04,Kel05,BO07}\ldots
These problems were tackled using mainly {\em topological} techniques.
The construction of a spectral triple is a proposition for bringing {\em geometry}
into play.

In this article, we revisit the construction of Pearson and Bellissard for the
transversals of some tiling spaces.
For this purpose, we use the formalism of Bratteli diagrams instead of Michon
trees.
This approach is equivalent, and applies in general to any ultrametric Cantor
set.
But in some cases, the diagram conveniently encodes the self-similarity, and is
very well suited to handle explicit computations.
For this reason, we will focus on diagrams arising from substitution tilings,
a class of tilings which is now quite well studied
\cite{Que87,Ken96,Sol97,AP98}.

Bratteli diagrams were introduced in the seventies for the classification of
$AF$-algebras \cite{Bra71}.
They were used by Ver\v{s}ik to encode measurable $\ZM$-actions, as a tool to
approach some dynamical systems by a sequence of periodic dynamical systems
\cite{VL92}.
They were then adapted in the topological setting to encode $\ZM$-actions on
the Cantor set \cite{HPS92,GPS95,DHS99}.
Then, Bratteli diagrams were used to represent the orbit equivalence relation
arising from an action of $\ZM^2$ \cite{GMPSa08} (and recently $\ZM^d$
\cite{GMPSb08}) on a Cantor set.
Yet, it is not well understood how the \emph{dynamics} itself should be
represented on the diagram.
The case we will look at is when the Cantor set is the transversal of a tiling
space, and the action is related to the translations.
The idea of parametrizing tilings combinatorialy dates back to the work of Gr\"{u}nbaum and Shephard in the seventies, and the picture of a Bratteli diagram can be
found explicitely in the book of Connes \cite{Con90} for the Penrose tiling.
However, it took time to generalize these ideas, and to understand the
topological and dynamical underlying questions.



\begin{figure}[!h]
\[
\xymatrix{
&   &  \bullet \ar@{-}[rrdd] & &  \bullet \ar@{-}[rrdd] & &  \bullet \ar@{.}[dr] & \\
&\circ \ar@{-}[ur] \ar@{-}[dr] & & & & & & \\
&   &  \bullet \ar@{-}[rr]  \ar@{-}[uurr] & &  \bullet\ar@{-}[rr] \ar@{-}[uurr] & &  \bullet \ar@{.}[r] \ar@{.}[ur] & \\
}
\]
\caption[A self-similar Bratteli diagram]%
{{\small A self-similar Bratteli diagram associated with the matrix
$\left(
\begin{matrix}
1 & 1 \\
1 & 0
\end{matrix}
\right)$ (root on the left).}}
\label{tLap09.fig-BdiagFibo}
\end{figure}

In the self-similar case (for example when the Cantor set is the
transversal of a substitution tiling space), the diagram only depends
on an adjacency (or Abelianization) matrix.
There is a natural \Cs associated with this matrix, called a
Cuntz--Krieger algebra \cite{CK80}. Its generators implement recursion
relations and therefore provide a method of computation for the eigenvalues of
the Laplace--Beltrami operators.

While Bratteli diagrams are suited to facilitate computations
for any self-similar Cantor set, we focus on diagrams associated with
substitution tilings.
Indeed, a tiling space comes with a convenient distance, which is encoded (up
to Lipschitz-equivalence) in a natural way by weights on the diagram.

\subsection*{Results of the Paper}

Let $\Bb$ be a weighted Bratteli diagram (Definition~\ref{tLap09.def-bratteli}, and~\ref{tLap09.def-weight}), and $(\partial \Bb,d_w)$ the ultrametric Cantor set
of infinite paths in $\Bb$ (Proposition~\ref{tLap09.prop-ultrametric}).

The Dixmier trace associated with the spectral triple gives a measure
$\mu_{\textrm Dix}$ on $\partial \Bb$.
The construction of Pearson--Bellissard gives a family of Laplace--Beltrami
operators $\Delta_s, \, s\in \RM$,  on \(L^2(\partial \Bb, d\mu_{\text{\rm
Dix}})\), see Definition \ref{tLap09.eq-LapExplicit}.
For all $s$, $\Delta_s$ is a nonpositive, self-adjoint, and unbounded operator.
For a  path $\gamma$, we denote by $[\gamma]$ the clopen set  of infinite paths with prefix $\gamma$, and by $\chi_\gamma$ its characteristic function.
 And we let $\ext(\gamma)$ denote the set of ordered pairs of edges that extend
 $\gamma$ one generation further.
 The operator $\Delta_s$ was shown to have pure point spectrum in \cite{PB09}.
 In this paper, we determine all its eigenelements explicitly.

 \begin{varthm}[Theorem~\ref{tLap09.thm-SpectrumLap}]
 The eigenspaces of $\Delta_s$ are given by the subspaces
 \[
  E_\gamma =
  \Bigl<
  \frac{1}{\mu_{\text{\rm Dix}}[\gamma\cdot a]} \chi_{\gamma\cdot a} -
  \frac{1}{\mu_{\text{\rm Dix}}[\gamma\cdot b]} \chi_{\gamma\cdot b} \, : \,
  (a,b) \in \ext(\gamma)
  \Bigr>
 \]
 for any finite path $\gamma$ in $\Bb$.
 We have  \( \dim E_\gamma = n_\gamma -1\) where $n_\gamma$ is the number of
 edges in $\Bb$ extending $\gamma$ one generation further.
 \end{varthm}

\begin{figure}[!h] 
 \[
  \xymatrix{
  & \circ \ar@{.}[dl] \ar@{-}[dr] &                         
  &&& & \circ \ar@{.}[dl] \ar@{-}[dr] &  \\
  \bullet \ar@{.}[dd] \ar@{.}[ddrr] & & \bullet \ar@{-}[ddll]^(.7){\gamma}       
  &&&  \bullet \ar@{.}[dd] \ar@{.}[ddrr] & & \bullet \ar@{-}[ddll]  \\
  & & &&& & & \\
  \bullet \ar@{.}[dd] \ar@{.}[ddrr] & & \bullet \ar@{.}[ddll]               
  &&& \bullet \ar@{-}[dd]_(.7){\gamma \cdot a} \ar@{-}[ddrr]^(.7){\gamma \cdot b} & & \bullet \ar@{.}[ddll] \\
  & & &&& & & \\
  & & &&& & &
  }
  \]
 \caption[haha]{{\small Example of eigenvectors of $\Delta_s$ for the Fibonacci
 diagram.}}
 \label{tLap09.fig-EfnFibo}
 \end{figure}
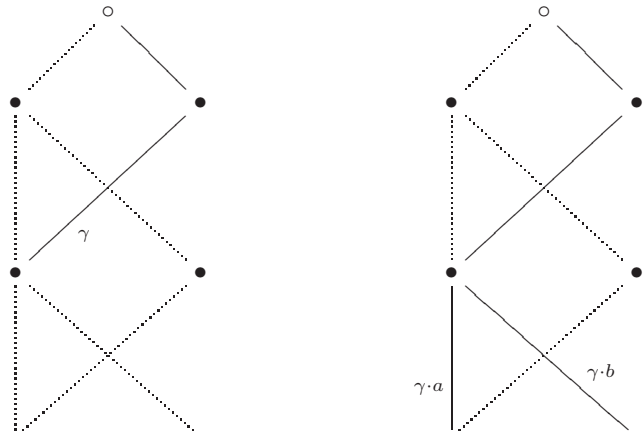

 The associated eigenvalues $\lambda_\gamma$ are also calculated explicitly, see
 equation \eqref{tLap09.eq-eigenfn}.
 An eigenvector of $\Delta_s$ is simply a weighted sum of the characteristic
 functions of two paths of the same lengths that agree apart from their last
 edge, see Figure \ref{tLap09.fig-EfnFibo} for an example.

 The transversal $\Xi$ to a substitution tiling space $(\Omega, \omega)$ of
 $\RM^d$, can be described by a stationary Bratteli diagram, like the Fibonacci
 diagram shown in Figure~\ref{tLap09.fig-BdiagFibo}.
 There is a ``natural map'', called the Robinson map in \cite{Kel95}, \( \psi :
 \partial \Bb \rightarrow \Xi\), which under some technical conditions
 (primitivity, recognizability, and border forcing) is a  homeomorphism
 (Theorem~\ref{tLap09.prop-homeo}).
 We endow $\Xi$ with  the {\em combinatorial metric} $d_\Xi$: two tilings are
 $\eps$-close if they agree on a ball of radius $1/\eps$ around the origin.
 Let $A$ be the Abelianization matrix of the substitution, and $\pf$ its
 Perron--Frobenius eigenvalue.
We denote by $w(\gamma)$ the weight of a finite path $\gamma$ in $\Bb$ (Definition~\ref{tLap09.def-weight}).

 \begin{varthm}[Theorem~\ref{tLap09.thm-bilip}]
 If there are constants $c_+>c_->0$ such that 
\( c_- \pf^{-n/d} \le w(\gamma) \le c_+ \pf^{-n/d}\) 
for all paths $\gamma$ of lengths $n$,
then the homeomorphism \(\psi:\partial\Bb\rightarrow\Xi\) is
($d_w$--$d_\Xi$) bi-Lipschitz.
\end{varthm}

In our case of substitution tiling spaces, we have a $\zeta$-function which is
given as in \cite{PB09} by
\[
\zeta(s) = \sum_{\gamma \in \Pi} w(\gamma)^s\,,
\]
where $\Pi$ is the set of finite paths in $\Bb$.
It is proven in \cite{PB09} that, when it exists, the abscissa of convergence
$s_0$ of $\zeta$ is the upper box dimension of the Cantor set.
For self-similar Cantor sets it always exists and is finite.

\begin{varthm}[Theorem~\ref{tLap09.thm-zeta}]
For a weighted Bratteli diagram associated with a substitution tiling space of
dimension $d$, the abscissa of convergence of the $\zeta$-function is $s_0=d$.
\end{varthm}

We also have an interpretation of $s_0$ which is not topological.
We link $s_0$ to the exponent of the complexity function.
This function $p$, associated with a tiling, counts the number of
distinct patches: $p(n)$ is the number of patches of radius $n$ (up to
translation).
We present two results.

\begin{varthm}[Theorem~\ref{tLap09.theo-box-cplex}]
For the transversal of \emph{any} minimal aperiodic tiling space with
a well-defined complexity function, the box dimension, when it exists,
is given by the following
limit:
\[
\dim \Xi = \limsup_{n\rightarrow +\infty}{\frac{\ln(p(n))}{\ln(n)}} \, .
\]
\end{varthm}
\noindent And we deduce the following.

\begin{varthm}[Corollary~\ref{tLap09.thm-complex}]
Let $\Xi$ be the transversal of a substitution tiling of dimension $d$, with
complexity function $p$. Then there exists a function $\nu$ such that:
\[
p(n) = n^{\nu(n)}, \qquad \textrm{with } \lim_{n\rightarrow
+\infty}{\nu(n)}=d\, .
\]
\end{varthm}

With the above choice of weights we can compute the Dixmier trace $\mu_{\textrm Dix}$.
Furthermore, there is a uniquely ergodic measure on $(\Omega,\RM^d)$, which was
first described by Solomyak~\cite{Sol97}.
It restricts to a measure $\mu_\Xi^{\textrm erg}$ on the transversal, and we have the following.

\begin{varthm}[Theorem~\ref{tLap09.thm-measure}]
With the above weights one has
\( \psi_\ast ( \mu_{\textrm Dix}) \  = \ \mu_\Xi^{\textrm erg} \).
\end{varthm}

As $\Bb$ is stationary, the sets of edges between two generations (excluding the
root) are isomorphic. 
Let us denote by $\Emod$ this set, and by $\Ee_0$  the set of edges linking to
the root.
Thanks to the self-similar structure we can define affine maps $u_e$,
\( e\in \Emod\), that act on the eigenvalues of $\Delta_s$ as follows (see section
\ref{tLap09.ssect-CKalg})
\[
u_e(\lambda_\eta) = \lambda_{U_e \eta} =
             \pf^{(d+2 - s)/d} \lambda_\eta + \beta_e \,,
\]
where $U_e \eta$ is an extension of the path $\eta$ (see Definition~%
\ref{tLap09.eq-actionCK-l2}),
and $\beta_e$ a constant that only depends on $e$.
Those maps $U_e, e\in \Emod$, form a representation of a Cuntz--Krieger algebra
associated with the matrix $A$.
If $\gamma = (\eps, e_1, e_2, \cdots e_n)$, $\eps \in \Ee_0$, $e_i\in \Emod$, is a
path of length $n$,
let $u_\gamma$ be the map \(u_{e_1} \circ u_{e_2} \circ \cdots u_{e_n}\).
Let us denote by $\lambda_\eps, \eps \in \Ee_0$, the eigenvalues of $\Delta_s$
corresponding to paths of length $1$.
For any other eigenvalue $\lambda_\gamma$, there is a (unique) $\lambda_\eps$
such that
\[
 \lambda_\gamma = u_\gamma (\lambda_\eps) =
\, \Lambda_s^n \lambda_\eps \ + \  \sum_{j=1}^n \, \Lambda_s^{j-1} \,\beta_{e_j} \,,
\]
where \(\Lambda_s = \pf^{(d+2-s)/d}\).
That is, the Cuntz--Krieger algebra allows to calculate explicitly the full
spectrum of $\Delta_s$ from the finite data of the \(\lambda_\eps, \eps \in
\Ee_0\), and \( \beta_e, e\in \Emod\) ---~which are immediate to compute, see
Section~\ref{tLap09.ssect-CKalg}.

\begin{figure}[htbp]
  \begin{center}
  \input{repart5.pstex_t}
  \end{center}
\caption{\small Distribution of the eigenvalues for the Fibonacci diagram.}
\label{tLap09.fig-fibo}
\end{figure}
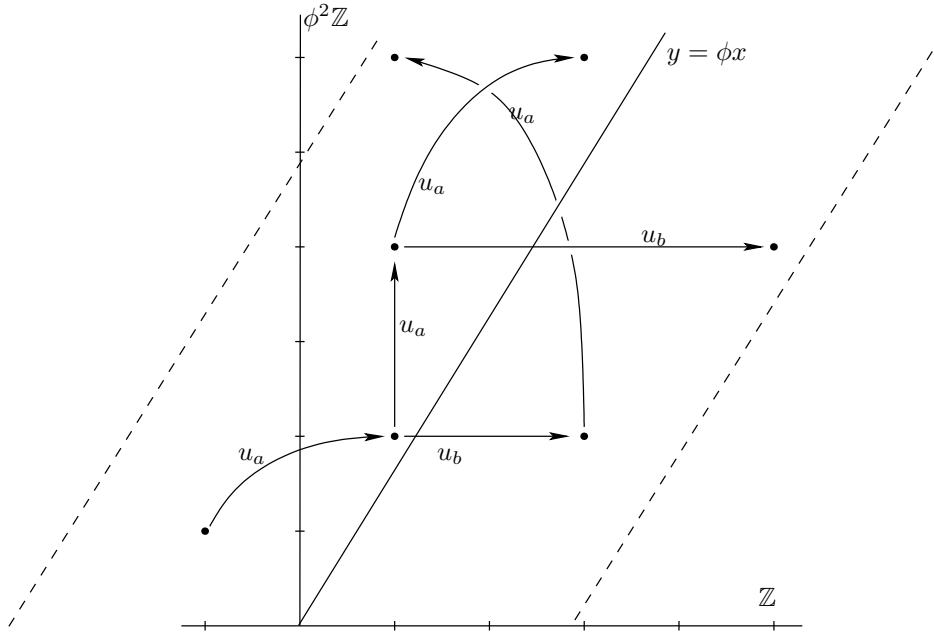

For instance, for the Fibonacci diagram (Figure \ref{tLap09.fig-BdiagFibo}) and $s=s_0=d$, there are only two such
$u_e$ maps, namely \(u_a(x) =  x \cdot \phi^2 - \phi\), and \(u_b(x) =  x \cdot
\phi^2 + \phi\), where \(\phi=(1+\sqrt{5})/2\) is the golden mean.
The eigenvalues of $\Delta_d$ are all of the form $ p + q \phi^2$ for integers
$p, q \in \ZM$.
They can be represented as points $(p,q)$ in the plane; these points stay within
a bounded distance to the line directed by the Perron--Frobenius eigenvector of the Abelianization matrix, see Figure~%
\ref{tLap09.fig-fibo}.
This is an example of a general result, valid for hyperbolic substitutions, see
Theorem~\ref{tLap09.thm-dist-eval}.
We treat further examples in Section~\ref{tLap09.sect-Ex}, in particular the
Thue--Morse and Penrose tilings.

The name Laplace--Beltrami operator for $\Delta_s$ can be justified by the
following two results.
For $s=s_0=d$, in analogy with the Laplacian on a compact $d$-manifold,
$\Delta_d$ satisfies the classical Weyl asymptotics, and the trace of its heat
kernel follows the leading term of the classical Seeley expansion.

Let $\Nn(\lambda)$ be the number of eigenvalues of $\Delta_d$ of modulus less
than $\lambda$.

\begin{varthm}[Theorem~\ref{tLap09.thm-Weyl}]
There are constants $0<c_-<c_+$ such that as \(\lambda \rightarrow + \infty\)
one has
\[
 c_- \lambda^{d/2} \le \Nn(\lambda) \le  c_+ \lambda^{d/2}\,.
\]
\end{varthm}

\begin{varthm}[Theorem \ref{tLap09.thm-Seeley}]
There are constants $0<c_-<c_+$ such that as \(t \downarrow 0\) one has
\[
c_- t^{-d/2} \le \TR \bigl( e^{t\Delta_d} \bigr) \le c_+ t^{-d/2}\,.
\]
\end{varthm}

\begin{acks}
The authors were funded by the NSF grants no.~DMS-0300398 and no.~DMS-0600956 of
Jean Bellissard, while visiting the Georgia Institute of Technology for the
Spring term 2009.
It is a pleasure to thank J.~Bellissard for his many supportive, insightful,
and enthusiastic remarks on this work.
The authors would like to thank Ian Putnam for useful discussions (which lead
in particular to the study of the eigenvalue distribution in section
\ref{tLap09.ssect-Evdist}), as well as Johannes Kellendonk and Marcy Barge for
inviting them to Montana State University, and useful discussions they had with
them there.
J.S. also aknowledges financial support from the SFB 701, Universit\"at Bielefeld, and would like to thank Michael Baake for his generous invitations.
\end{acks}

\section{Weighted Bratteli diagrams and substitutions}
\label{tLap09.sect-wBdiag}

We first give general the definitions of Cantor sets, Bratteli diagrams,
and substitution tilings. We then describe how to associate diagrams to
tilings.

\begin{defini}
A Cantor set is a compact, Hausdorff, metrizable topological space, which
is totally disconnected and has no isolated points.

An ultrametric $d$ on a topological space $X$ is a metric which satisfies this
strong triangle inequality:
\[
\forall x,y,z \in X, \quad d(x,y) \leq \max\{d(x,z),d(y,z)\}.
\]
\end{defini}


\subsection{Bratteli diagrams}

To a substitutive system, one can naturally associate a combinatorial object
named a Bratteli diagram.
These diagrams were first used in the theory of $C^*$-algebras, to classify
AF-algebras.
Then, it was mainly used to encode the dynamics of a minimal action of $\ZM$
on a Cantor set.

\begin{defini}\label{tLap09.def-bratteli}
A Bratteli diagram is an oriented graph defined as follows:
\[
\Bb = (\Vv, \Etot, r,s),
\]
Where $\Vv$ is the set of vertices, $\Etot$ is the set of directed edges,
and $r,s$ are functions $\Etot \rightarrow \Vv$ (range and source),
which define adjacency.
We have a partition of $\Vv$ and $\Etot$ in finite sets:
\[
\Vv = \bigcup_{n \geq 0}{\Vv_n} \quad ; \quad \Etot = \bigcup_{n \geq 0}{\Ee_n},
\]
$\Vv_0$ is a single element called the root and noted $\circ$.
The edges of $\Ee_n$ have their source in $\Vv_n$ and range
in $\Vv_{n+1}$, that is:
\[
r: \Ee_n \rightarrow \Vv_{n+1}, \quad s:\Ee_n \rightarrow \Vv_n \, .
\]
We ask that $s^{-1}(v) \neq \emptyset$ for all $v \in \Vv$, and
$r^{-1}(v) \neq \emptyset$ for all $v \in \Vv\setminus\Vv_0$.
\end{defini}

\begin{rem}
\label{tLap09.rem-tree}
If, for all $v \in \Vv \setminus\{\circ\}$, $r^{-1} (v)$ is a single
edge, then $\Bb$ is a tree.
In that sense, the formalism of Bratteli diagrams includes the case of trees,
and so is a generalization of the case studied in \cite{PB09}.
However, our goal is to restrict to self-similar diagrams, for which
computations are easier.
\end{rem}

\begin{defini}
A path $\gamma$ of length $n \in \NM \cup \{+\infty\}$ in a Bratteli diagram
is an element
\[
(\eps_0, e_1, e_2,\ldots) \in \prod_{i=0}^{n-1}{\Ee_i},
\]
which satisfies:
\[
\text{for all } 0 \leq i < n, \quad r(e_i) = s(e_{i+1}).
\]
We call $\Pi_n$ the set of paths of length $n < +\infty$, $\Pi$ the set of all
finite paths, and $\partial \Bb$ the set of infinite paths.
\end{defini}
The function $r$ naturally extends to $\Pi$: if
$\gamma = (\eps_0, \ldots, e_n)$, then $r(\gamma) := r(e_n)$.

In addition to the definition above, we ask that a Bratteli diagram satisfies
the following condition.
\begin{hypo}
\label{tLap09.hypo-split}
For all $v \in \Vv$, there are at least two distinct infinite paths through
$v$.
\end{hypo}

\begin{rem}
\label{tLap09.rem-simple-edges}
If, for all $n \geq 1$ and $(v, v') \in \Vv_n \times \Vv_{n+1}$, there is at
most one edge from $v$ to $v'$, we can simply encode path by
vertices: the following map is an homeomorphism onto its image.
\[
(\eps_0, e_1, \ldots) \longmapsto (\eps_0, r(e_1), \ldots).
\]
\end{rem}

\begin{defini}
Given two finite or infinite paths $\gamma$ and $\gamma'$, we note
$\gamma \wedge \gamma'$ the (possibly empty) longest common prefix of $\gamma$
and $\gamma'$.
\end{defini}

The set $\partial \Bb$ is called the boundary of $\Bb$. It has a natural
topology inherited from the product topology on $\prod_{i=0}^{+\infty}{\Ee_i}$,
which makes it a compact and totally disconnected set.
A basis of neighborhoods is given by the following sets:
\[
[\gamma] = \{ x \in \partial \Bb \ ; \ \gamma \text{ is a prefix of } x\}.
\]

Hypothesis~\ref{tLap09.hypo-split} is the required condition to make sure that
there are no isolated points. This implies the following.
\begin{proposi}
With this topology, $\partial \Bb$ is a Cantor set.
\end{proposi}

\begin{defini}
\label{tLap09.def-weight}
A weight on $\partial \Bb$ is a function $w: \Vv \rightarrow \RM_+^*$,
satisfying the following conditions:
\begin{enumerate}[(i)]
\item $w(\circ) = 1$;
\item $\sup\{w(v) \ ; \ v \in \Vv_n\}$ tends to $0$ when $n$ tends to infinity;
\item $\forall e \in \Etot$, $w(s(e)) > w(r(e))$.
\end{enumerate}
\end{defini}
A weight extends naturally on paths: by definition, $w(\gamma) := w(r(\gamma))$.

\begin{proposi}
\label{tLap09.prop-ultrametric}
We define a function $d_w$ on $(\partial \Bb)^2$ by:
\[
d_w(x,y) = \left\{ \begin{array}{ll}
              w\big( r(x \wedge y) \big) & \text{if } x \neq y; \\
              0 & \text{otherwise}
           \end{array}\right.
\]
It is a ultrametric on $\partial \Bb$, which is compatible with the topology
defined above, and so $(\partial \Bb,d_w)$ is a ultrametric Cantor set.
\end{proposi}

A case of interest is when the Bratteli diagram is self-similar.
For a diagram, self-similarity means that
all the $\Vv_n$ are isomorphic, and all $\Ee_n$ are isomorphic for
$n \geq 1$, and $r,s$ commute with these identifications.
We will focus on self-similarity when the diagram is associated with
a substitution.

\subsection{Substitution tilings}
\label{tLap09.ssect-tilings}

Let us give the definition of a tiling of $\RM^d$.
The tilings we are interested in are constructed from a prototile set and an
inflation and substitution rule on the prototiles.
The notion of substitution dates back to the sixties.
For example the self-similarity of the Penrose tiling already appeared in \cite{Pen79}.
However, systematic formalisation of substitution tilings and properties of
their associated tiling space was done in \cite{Ken96,Sol98} in the nineties.
In this section, we follow the description by Anderson and Putnam~\cite{AP98}
with some minor changes.
In particular, some non trivial facts are cited along the text. We do not
claim to cite the original authors for all of these.
The reader can refer to the reviews \cite{Rob04} and \cite{Fra08}.

By \emph{tile}, we mean a compact subset of $\RM^d$, homeomorphic to the closed
unit ball.
A tile is \emph{punctured} if it has a distinguished point in its
interior.
A prototile is the equivalence class under translation of a tile.
Let $\Aa$ be a given set of punctured prototiles. All tilings will be made from
these tiles.

A patch $p$ on $\Aa$ is a finite set of tiles which have disjoint interior.
We call $\Aa^*$ the set of patches modulo translation.
A partial tiling is an infinite set of tiles with disjoint interior,
and the union of which is connected, and a tiling $T$ is a partial tiling which
covers $\RM^d$, which means:
\[
\Supp(T) := \bigcup_{t \in T}{t} = \RM^d \, .
\]

A substitution rule is a map $\omega$ which maps tiles to patches, and such that
for all tile $p \in \Aa$, $\Supp(\omega (p)) = \lambda \Supp (p)$ for some
$\lambda > 1$.
The factor $\lambda$ inflates the tiles $p$, which is then cut into pieces;
these pieces are elements of $\Aa$.
The map $\omega$ extends to patches, partial tilings and tilings.

Since we will represent tilings spaces by diagrams, we have a specific interest
for the combinatorics of the substitution.
\begin{defini}
Given a substitution $\omega$, its \emph{Abelianization matrix} is an
integer-valued matrix $A_\omega = (a_{pq})_{p,q \in \Aa}$ (or simply $A$)
defined by:
\[
\forall p,q \in \Aa, \quad a_{pq} = \textrm{ number of distinct translates
of $p$ included in $\omega(q)$.}
\]
\end{defini}


We now define the tiling space associated with $\omega$.
\begin{defini}
The tiling space $\Omega$ is the set of all tilings $T$ such that for all patch
$p$ of $T$, $p$ is also a subpatch of some $\omega^n(t)$ ($n \in \NM$, $t \in
\Aa$).
\end{defini}

We make the following assumptions on $\omega$:
\begin{hypo}
\begin{enumerate}[(i)]
\item The matrix $A$ is \emph{primitive}: $A^n$
has non-negative entries for some $n$. 
\item \emph{Finite local complexity} (FLC): for all $R > 0$, the set of all patches
$p \subset T$ which can be included in a ball of radius $R$ is finite up to
translation, for all $T \in \Omega$.
\end{enumerate}
\end{hypo}

Just as every tile has a distinguished point inside it, $\Omega$ has a
distinguished subset. This is the Cantor set which will be associated with
a Bratteli diagram.
\begin{defini}
Let $\Xi$ be the subset of $\Omega$ of all tilings $T$ such that $0$ is the
puncture of one of the tiles of $T$.
This set is called the canonical transversal of $\Omega$.
\end{defini}
Note that, by this definition, the substitution $\omega$ extends to a map
$\Omega \rightarrow \Omega$.

We will assume in the following that $\omega:\Omega \rightarrow \Omega$ is
one-to-one.
This condition is equivalent (Solomyak~\cite{Sol98}) to the fact that no tiling
in $\Omega$ has any period.
Furthermore, it implies that $\Omega$ is not empty, and $\omega: \Omega
\rightarrow \Omega$ is onto.


The set of interest for us is $\Xi$. Its topology is given as follows.
For any patch $p$, define the following subset of $\Xi$:
\[
U_p = \{T \in \Xi \ ; \ p \subset T \}.
\]
Note that $U_p$ can be empty.
Nevertheless, the family of all the $U_p$'s is a basis for a topology on $\Xi$.


Let us now define a distance $d$ on $\Omega$:
\begin{multline*}
d(T,T') = \inf\bigg(\Big\{ \epsilon > 0 \ ; \ \exists x,y \in B(0,\epsilon) \text{ such
that } \\ B(0,1/\epsilon) \subset \Supp\big( (T-x) \cap (T'-y) \big) \Big\}
\cup \{1\}\bigg).
\end{multline*}
Two tilings are $d$-close when, up to a small translation, they agree on a large
ball around the origin.
This distance, when restricted to $\Xi$, is compatible with the topology
defined above.
With this topology, $\Xi$ is a Cantor set.

Furthermore, the map $\omega: \Omega\rightarrow\Omega$ is a homeomorphism,
and the dynamical system $(\Omega,\RM^d)$ given by translation is continuous
and uniquely ergodic~\cite{Sol97}.

Substitution tiling spaces are \emph{minimal}, which means that every $\RM^d$-orbit
is dense in $\Omega$.
Combinatorially, this is equivalent to the fact that these tilings are repetitive:
for all $R > 0$, there is a bound $\rho_R$, such that every patch of size $R$ appears
in the tiling whithin range $\rho_R$ to any tile.

We need an additional assumption ---~the border forcing condition. It is
required in order to give a good representation of $\Xi$ as the boundary of a
Bratteli diagram.

\begin{defini}
\label{tLap09.def-unambiguous}
Let $p \in \Ll$ be a patch. The \emph{maximal unambiguous extention}
of $p$ is the following patch, obtained as an intersection of tilings:
\[
\Ext(p) = \bigcap\Big\{ T \in \Omega \ ; \ p \subset T \Big\}.
\]
\end{defini}
So any tiling which contains $p$ also contains $\Ext(p)$.
The patch $\Ext(p)$ is called the \emph{empire} of the patch $p$ by some
authors \cite{GS89}.

\begin{defini}
\label{tLap09.def-bordforc}
Assume there is some $m \in \NM$ and some $\epsilon > 0$ such that for all
$t \in \Aa$,
\[
\Ext\big(\omega^m(t)\big) \text{ contains an $\epsilon1$-neighborhood of }
\omega^m(t).
\]
Then $\omega$ is said to \emph{force its border}.
\end{defini}

It is always possible to give labels to the tiles, in order to change a
non-border forcing substitution into a border-forcing one. See~\cite{AP98}.

\begin{exam}
\label{tLap09.ex-fiboconj}
Let us give an example in dimension one. Consider the substitution defined on
symbols by:
\[
\begin{array}{ccc}
a & \mapsto & baa \\
b & \mapsto & ba
\end{array}
\]
It has a geometric realization, with $a$ being associated with an interval of
length $\phi = (1+\sqrt{5})/2$, and $b$ with an interval of length $1$.
Then the map above is a substitution in the sense of our definitions, with
$\lambda = \phi$.
The transversal $\Xi$ is naturally identified with a subset of $\{a,b\}^\ZM$
(the \emph{subshift} associated with $\omega$).
The map $\omega$ satisfies the border forcing condition.
Indeed, if $p$ is a patch of a tiling $T$
(which we identify symbolically with a finite word on $\{a,b\}$ in a
bi-infinite word), let $x,y \in \{a,b\}$ be the letters preceding
and following $p$ in $T$ respectively.
Then $\omega(x)$ ends by an $a$ and $\omega(y)$ begins by a $b$.
Therefore, $\omega(p)$ is always followed by $b$ and preceded by $a$.
This proves that $\omega$ forces its border.
\end{exam}

\subsection{Bratteli diagrams associated with substitution tilings}

We show how to associate a Bratteli diagram to a substitution,
and identify the transversal with the boundary of the diagram.
It is clear that these two are homeomorphic, being Cantor sets.
We will give an explicit and somehow natural homeomorphism.


Let $\omega$ be a primitive, FLC and border forcing
substitution, and $A = (a_{ij})_{i,j \in \Aa}$ its Abelianization matrix.
Let $\lambda$ be the expansion factor associated with $\omega$.
Let $\Xi$ be the transversal of the tiling space associated with $\omega$.

The diagram $\Bb$ associated with the substitution is defined as follows.
\begin{defini}
\label{tLap09.def-brat-subst}
The diagram associated with $\omega$ is the diagram $\Bb = (\Vv,\Etot,r,s)$,
where:
\[
\forall n \geq 1, \quad \Vv_n = \Aa \times \{n\},
\]
for all $n \geq 1$, for all $i,j \in \Aa$, there are exactly
$a_{ij}$ edges of $\Ee_n$ from $(i,n)$ to $(j,n+1)$,
and for all $v \in \Vv_1$ there is exactly one edge $\eps_v$ from the root
$\circ$ to $v$.
\end{defini}

We call $\Emod$ the set of ``models'' of edges from one generation to another.
All $\Ee_n$ are copies of $\Emod$ (for example with identification
$\Ee_n = \Emod \times \{n\}$).
Since the models for the vertices are the elements of $\Aa$, we do not use a
specific name.
The maps source and range can be restricted as maps $\Emod \rightarrow \Aa$.
When the ``depth'' of a vertex is not important,  we will sometimes consider
$r,s$ as functions valued in $\Aa$.

\begin{rem}
Combinatorially, the diagram only depends on the Abelianization
matrix of $\omega$.
It is indeed possible to associate a diagram to a primitive matrix $A$ with
integer coefficients.
We would have similar definitions, with $\Vv_n = I \times \{n\}$, with
$I$ the index set of the matrix.
\end{rem}

\begin{proposi}
This is a Bratteli diagram in the sense of Definition~\ref{tLap09.def-bratteli},
except in the case $A = (1)$.
\end{proposi}


\begin{exam}
\label{tLap09.rem-fibo}
Figure~\ref{tLap09.fig-BdiagFibo} is an example of the self-similar diagram
associated with the Fibonacci substitution: $a \mapsto ab$; $b \mapsto a$.
\end{exam}

There is a correspondence between the paths on $\Bb$ and
the transversal $\Xi$.
It depends on a choice on the edges which remembers the geometry of the
substitution: each $e \in \Emod$ from $a \in \Aa$ to $b \in \Aa$ corresponds to a
different occurrence of $a$ in $\omega(b)$.
This correspondence is a homeomorphism $\phi: \Xi \rightarrow \partial \Bb$,
called the Robinson map, as defined in~\cite{Kel95}.
We first give the definition of $\phi$, and we will then give a condition the
weight function $w$ so that $\phi$ is bi-Lipschitz.

To construct $\phi$, start with $T \in \Xi$.
Let $t \in T$ be the tile containing $0$, and $[t] \in \Aa$ the corresponding
prototile.
Then, the first edge of $\phi (p)$ is the edge from $\circ$ to
$([t],1) \in \Vv_1$.
Assume that the prefix of length $n$ of $\phi(T)$ is already constructed and
ends at vertex $([t'],n)$,
where $t'$ is the tile of $\omega^{-(n-1)}(T)$ which contains $0$.
Let $t''$ be the tile containing $0$ in $\omega^{-n}(T)$;
then the $(n+1)$-th edge of $\phi(T)$ is the edge corresponding to
the inclusion of $t'$ in $\omega(t'')$. This edge ends at $([t''],n+1)$.
By induction we construct $\phi(T)$.

\begin{theo}[Theorem~4 in \cite{Kel95}]
\label{tLap09.prop-homeo}
The function $\phi : \Xi \rightarrow \partial \Bb$ is an homeomorphism.
\end{theo}

We give explicitly the inverse for $\phi$ as follows.
Let us define $\psi: \Pi \rightarrow \Aa^*$ which is increasing in the sense
that if a path is a prefix of another, then the patch associated to the first
is included in the patch associated to the second.
The image of an infinite path will then be defined as the union of the images
of its prefixes.

If $(\eps)$ is a path of length one, then $\psi(\eps)$ is defined as the tile
$r(\eps)$ with puncture at the origin.
Now, given a path $\gamma$ of length $n$, assume its image by $\psi$ is well
defined, and is some translate of $\omega^{n-1}(r(\gamma))$, with the origin
at the puncture of one of its tiles.
Consider the path $\gamma.e$, with $e \in \Ee_n$ such that $s(e) = r(\gamma)$.
Then, $e$ encodes an inclusion of $s(e)$ inside $\omega(r(e))$.
It means that  it encodes an inclusion of $\omega^{n-1}(r(\gamma))$ inside
$\omega^n(r(e))$.
The patch $\psi(\gamma .e)$ is defined as the translate of $\omega^n(r(e))$
such that the inclusion $\psi (\gamma) \subset \psi (\gamma.e)$
is the inclusion defined by the edge $e$.

Remark that if $x \in \partial \Bb$, $\psi(x)$ does not \emph{a priori}
define more than a partial tiling.
The fact that it corresponds to a unique tiling results from the border forcing
condition.
Then the map $\psi$, extended to infinite paths, is the inverse of $\phi$.


Continuity can be proved directly, but can also be seen as a consequence of
Theorem~\ref{tLap09.thm-bilip}, which we prove later.

\begin{lemma}
\label{tLap09.lemma-border}
There exists $C_1,C_2$ positive constants, such that
for all path $\gamma$ of length $n$ in $\Bb$,
\begin{eqnarray*}
B(0,C_1 \lambda^n) \subset \Supp\Big(\Ext\big(\psi (\gamma)\big)\Big), \\
B(0,C_2 \lambda^n) \not\subset \Supp\Big(\Ext\big(\psi (\gamma)\big)\Big),
\end{eqnarray*}
where $\lambda$ is the expansion factor of $\omega$.
\end{lemma}

\begin{proof}
Let $k$ be the smallest number such that for all $v \in \Vv$, there exists
two distinct paths of length $k$ starting from $v$.
For example, $k=1$ when there are two elements of $\Etot$ starting from $v$,
for all $v$, and $k=2$ in the case of Fibonacci,
pictured in Figure~\ref{tLap09.fig-BdiagFibo}.
Let $C$ be the maximum diameter of the tiles.
Let $\gamma$ be a path of length $n$.
Then, one can find two distinct extensions of length $n+k$ of $\gamma$;
call them $\gamma_1$ and $\gamma_2$.
Remember that $\psi(\gamma_i)$ is some  translate of $\omega^{n+k-1}(t_i)$ with
$t_i \in \Aa$, and $0 \in \Supp\big(\psi(\gamma)\big)$.
So $\psi(\gamma_1)$ and $\psi(\gamma_2)$ differ within range $C\lambda^{n+k-1}$.
Therefore, with $C_1 = C\lambda^{k-1}$, we have the second line.

For the first inclusion, let $m$ be the exponent for which $\omega$ satisfies
the border forcing condition, and $\gamma$ a path of length $m$.
Let $t \in \Aa$ be the range of $\gamma$.
Then $\Ext(\omega^m(t))$ covers an $\epsilon$-neighborhood of $\Supp(\omega^m(t))$.
It means that, since $\psi (\gamma)$ is a translate of $\omega^m(t)$ which contains $0$,
\[
B(0,\epsilon) \subset \Ext(\psi (\gamma)).
\]
Similarly, for all path of length $m+k$,
\[
B(0,\epsilon \lambda^k) \subset \Ext(\psi(\gamma)).
\]
With $n=m+k$ and $C_1 = \epsilon/\lambda^m$, one has the result.
For $n < m$, the inequality still holds (up to a reduction of $\epsilon$).
\end{proof}


\begin{rem}
\label{tLap09.rem-group}
We can  simplify the diagram in the case where a symmetry group
$G \subset \mathrm{O}_n(\RM^d)$ acts freely on $\Aa$, such that
\[
\forall g \in G, \ \forall p \in \Aa, \quad
\omega(g \cdot p) = g \cdot \omega(p).
\]
With this, it is possible to extend the action of $G$ to the edges of $\Etot$.
It induces naturally an action of $G$ on $\partial \Bb$.
The group $G$ also acts on $\Omega$ by isometries, and these two actions are
conjugate by $\psi$.

One can ``fold'' the diagram by taking a quotient
as follows.
Define:
\(
\Bb' = (\Vv', \Etot', r',s'),
\)
where $\Vv'_0 = \Vv$, $\Ee'_0 = \Ee_0 \times G$, and
\[
\forall n \geq 1, \quad \Vv'_n = \Vv_n / G \text{ and }
  \Ee'_n = \Ee_n / G.
\]
Furthermore, for $[e]_G \in \Ee_n/G$,
\(
r'([e]_G) := [r(e)]_G\) and \(s'([e]_G) := [s(e)]_G.
\)
These definitions do not depend on the choice of the representant $e$.
For $(e,g) \in \Ee'_0=\Ee_0 \times G$, define
\(
r'((e,g)) := [r(e)]_G\) and \(s'((e,g)) := \circ.
\)
One can check that these two diagrams are ``the same'' in the sense that
the tree structure of their respective sets of paths of finite length is the
same.

In terms of substitution, the image by $\psi$ of a path
$\gamma = (\eps_0, e_1, \ldots, e_{n-1})$ is the image of
$\omega^n (r(\gamma))$ under some element of $\RM^d \rtimes G$ .
The translation part is given by the truncated path $(e_1, \ldots, e_{n-1})$.
The rotation part is encoded in the first edge $\eps_0 \in \Ee'_0$:
$\eps_0 = (\eps',g)$, where $\eps'$ brings no additional information, but $g$
corresponds to the choice of an orientation for the patch.

There is still an action of $G$ on $\partial \Bb'$, defined by:
\(
g \cdot ((\eps_0,h), e_1, \ldots) = ((\eps_0,gh), e_1, \ldots)\, .
\)
\end{rem}


\subsection{Weights and metric}
\label{tLap09.ssect-weights}

We gave in Theorem~\ref{tLap09.prop-homeo} an explicit homeomorphism $\phi$
between the boundary of the Bratteli diagram, and the transversal of the tiling
space.
Since we are interested in metric properties of the Cantor set, we now show that
a correct choice of weights on the vertices of the Bratteli
diagram gives a metric on $\partial \Bb$ which is Lipschitz equivalent
to the usual metric on $\Xi$ (the equivalence being induced by $\phi$).

\begin{theo}
\label{tLap09.thm-bilip}
Let $\Bb$ the Bratteli diagram associated with a substitution $\omega$.
Let $\lambda$ be the inflation factor of $\omega$.
We make the following assumption on the weight function $w$:
\begin{equation}
\label{tLap09.eq-subs-weights}
\forall n \geq 1, \ \forall v \in \Vv, \quad w(v,n+1) = \frac{1}{\lambda}w(v,n).
\end{equation}
Then the function $\phi: (\Xi,d)\rightarrow (\partial \Bb,d_w)$ defined
in Proposition~\ref{tLap09.prop-homeo} is a bi-Lipschitz homeomorphism.
\end{theo}

\begin{proof}
It is enough to show that there are two constants $m,M>0$,
such that for all $\gamma \in \Pi$,
\[
m \leq \frac{\diam(U_\gamma)}{\diam(\phi^{-1}(U_\gamma))} \leq M.
\]
Since the $U_\gamma$ are a basis for the topology of $\partial\Bb$,
this will prove the result.

Let $\gamma \in \Pi_n$.
By Lemma~\ref{tLap09.lemma-border}, any two tilings in $\phi^{-1}(U_\gamma)$
coincide on a ball of radius at least $C_1 \lambda^n$.
Therefore,
\[
\diam(\phi^{-1}(U_\gamma)) \leq \frac{1}{C_1}\lambda^{-n}.
\]
On the other hand, it is possible to find two tilings in $\phi^{-1}(U_\gamma)$,
which disagree on a ball of radius $C_2 \lambda^n$.
Therefore,
\[
\diam(\phi^{-1}(U_\gamma)) \geq \frac{1}{C_2}\lambda^{-n}.
\]
And by definition of the weights,
\[
\min\{w(v,1) \ ; \ v \in \Vv \} \lambda^{-n+1}
\leq \diam (U_\gamma) \leq \max\{w(v,1) \ ; \ v \in \Vv \} \lambda^{-n+1}.
\]
Together with the previous two inequalities, this proves the result.
\end{proof}



\section{Spectral triple, $\zeta$-function, and complexity}
\label{tLap09.sect-ST}


\subsection{Spectral triple}
\label{tLap09.ssect-ST}

Let $\Bb$ be a weighted Bratteli diagram, and let $(\partial\Bb,d)$ be the
ultrametric Cantor set of infinite rooted paths in $\Bb$.
Pearson and Bellissard built in \cite{PB09} a {\em spectral triple} for
$(\partial\Bb,d)$ when $\Bb$ is a {\em tree} (that is when for all vertex
$v\in \Vv \setminus \{ \circ \}$, the fiber $r^{-1}(v)$ contains a single point).
In our setting, their construction is adapted as follows.

A {\em choice function} on $\Bb$ is a map
\begin{equation}
\label{tLap09.eq-extension}
\tau \ : \
\left\{
\begin{array}{ccc}
\Pi & \rightarrow & \partial \Bb \times \partial \Bb \\
\gamma & \mapsto & \bigl( \tau_+(\gamma), \tau_-(\gamma) \bigr)
\end{array}
\right.
\quad \text{\rm such that} \quad
d \bigl(
\tau_+(\gamma), \tau_-(\gamma) 
\bigr) = \diam [\gamma] \,,
\end{equation}
and we denote by $E$ the set of choice functions on $\Bb$.
Let \(C_{\text{\rm Lip}}(\partial \Bb)\) be the \Cs of Lipschitz continuous functions on $(\partial\Bb,d)$. 
Given a choice $\tau \in E$ we define a faithfull $\ast$-{\em representation} $\pi_\tau$ of $C_{\text{\rm Lip}}(\partial \Bb)$ by bounded operators on the Hilbert space \(\Hh = l^2(\Pi) \otimes \CM^2 \) as
\begin{equation}
\label{tLap09.eq-representation}
\pi_\tau (f) = \bigoplus_{\gamma \in \Pi}
\left[
\begin{array}{cc}
f\bigl(\tau_+(\gamma)\bigr) & 0 \\
0 & f\bigl(\tau_-(\gamma)\bigr)
\end{array}
\right] \,.
\end{equation}
This notation means that for all $\xi \in \Hh$ and all $\gamma \in \Pi$,
\[
\big(\pi_\tau (f) \cdot \xi\big)(\gamma) = \left[ \begin{array}{cc}
f\bigl(\tau_+(\gamma)\bigr) & 0 \\
0 & f\bigl(\tau_-(\gamma)\bigr)
\end{array}
\right]
\cdot \xi(\gamma)\, .
\]
A {\em Dirac operator} $D$ on $\Hh$ is given by
\begin{equation}
\label{tLap09.eq-Dirac}
D = \bigoplus_{\gamma \in \Pi} \frac{1}{\diam[\gamma]}
\left[
\begin{array}{cc}
0 & 1 \\
1 & 0 
\end{array}
\right] \,,
\end{equation}
that is $D$ is a self--adjoint unbounded operator such that \((D^2+1)^{-1}\) is compact, and the commutator
\begin{equation}
\label{tLap09.eq-commutator}
\bigl[
D, \pi_\tau (f)
\bigr]
= \bigoplus_{\gamma \in \Pi} 
\frac{f\bigl(\tau_+(\gamma)\bigr) - f\bigl(\tau_-(\gamma)\bigr)}{\diam[\gamma]}
\left[
\begin{array}{cc}
0 & -1 \\
1 & 0 
\end{array}
\right] \,,
\end{equation}
is bounded for all $f\in C_{\text{\rm Lip}}(\partial \Bb)$.
Finally a {\em grading operator} is given by 
\(\Gamma = 1_{l^2(\Pi)} \otimes \left[
\begin{array}{cc}
1 & 0 \\
0 & -1 
\end{array}
\right]\),
and satisfies \(\Gamma^2 = \Gamma^\ast = \Gamma\), and commutes with $\pi_\tau$ and anticommutes with $D$.
The following is Proposition 8 in \cite{PB09}.
\begin{proposi} 
\label{tLap09.prop-ST}
\( \bigl(C_{\text{\rm Lip}}(\partial \Bb), \Hh, \pi_\tau,D,\Gamma\bigr) \) is an {\em even spectral triple} for all $\tau \in E$. 
\end{proposi}

In \cite{PB09} the set of choice functions $E$ is considered an analogue of a tangent bundle over $\partial \Bb$, so that the above commutator is interpreted as the directional derivative of $f$ along the choice $\tau$.
The metric on $\partial \Bb$ is then recovered from the spectral triple by using Connes formula, {\it i.e.} by taking the supremum over all directional derivatives
\begin{theo}[Thm.~1 in \cite{PB09}]
The following holds:
\begin{equation}
\label{tLap09.eq-Connesdist}
d(x,y) = \sup 
\bigl\{
|f(x)-f(y)| \ ; \ f\in C_{\text{\rm Lip}}(\partial \Bb)\,, \ 
\sup_{\tau \in E} \|[D, \pi_\tau(f)] \| \le 1
\bigr\} \,.
\end{equation}
\end{theo}

\begin{defini}
The {\em $\zeta$-function} associated with the spectral triple is given by:
\begin{equation}
\label{tLap09.eq-zeta}
\zeta(s) = \frac{1}{2} \TR \bigl( |D|^{-s} \bigr) 
= \sum_{\gamma \in \Pi} \diam[\gamma]^s \,,
\end{equation}
\noindent and we will denote by $s_0 \in \RMbar$ its {\em abscissa of convergence}, when it
exists.
\end{defini}
 
We now assume that the weight system on $\Bb$ is such that $s_0 \in \RM$.

\begin{defini}
\label{tLap09.def-Dixmier}
The {\em Dixmier trace} of a function $f\in C_{\text{\rm Lip}}(\partial \Bb)$ is
given by the following limit, when it exists:
\begin{equation}
\label{tLap09.eq-Dixmier}
\mu(f) = \lim_{s\downarrow s_0} 
\frac{ \TR \bigl(|D|^{-s} \pi_\tau(f) \bigr) }{\TR \bigl(|D|^{-s}\bigr)} \, .
\end{equation}
\end{defini}

It defines a probability measure on $\partial \Bb$ and does not depend on the
choice $\tau \in E$ (Theorem~3 in \cite{PB09}).
Furthermore, when $f = \chi_\gamma$ is a characteristic function, the limit above can
be rewriten:
\begin{equation}
\mu([\gamma]) := \mu(\chi_\gamma) = \lim_{s \downarrow s_0}{
\frac{\sum_{\eta \in \Pi_\gamma}{w(r(\eta))^s}}{\sum_{\eta \in \Pi}{w(r(\eta))^s}} \, .
}
\end{equation}
where $\Pi_\gamma$ stands for the set of paths in $\Pi$ with prefix $\gamma$.

\subsection{The measure on $\partial \Bb$}
\label{tLap09.ssect-zetameasure}

Let $(\Omega, \omega)$ be a substitution tiling space, $A$ be the Abelianization
matrix of $\omega$, and $\Bb$ be the associated Bratteli diagram, as in
previous section (Definition~\ref{tLap09.def-brat-subst}).
We assume that $\Bb$ comes together with a weight function $w$, which satisfies
the properties given in Definition~\ref{tLap09.def-weight}, and adapted to
the substitution $\omega$.
In particular, it satisfies $w(a,n+1) = \lambda^{-1} w(a,n)$ for all $a \in \Aa$
and $n \in \NM$.

As we assumed $A$ to be primitive, it has a so called \emph{Perron--Frobenius}
eigenvalue, denoted $\pf$, which satisfies the following (see for
example~\cite{Hor94}):
\begin{enumerate}[(i)]
\item $\pf$ is strictly greater than $0$, and equals the spectral radius of $A$;
\item For all other eigenvalue $\nu$ of $A$, $\left\vert \nu \right\vert < \pf$;
\item The right and left eigenvectors, $v_R$ and $v_L$, have strictly positive
coordinates;
\item If $v_R$ and $v_L$ are normalized so that $\sca{v_R}{v_L} = 1$, then:
\[
\lim_{n \rightarrow +\infty}{\frac{A^n}{\pf^n}} = v_L v_R^{\tra};
\]
\item Any eigenvector of $A$ with non-negative coordinates corresponds to the
eigenvalue $\pf$.
\end{enumerate}

A classical result about linear dynamical systems together with the properties
above gives the following result, which will be needed later.
\begin{lemma}\label{tLap09.lem-rec}
Let $A$ be a primitive matrix with Perron--Frobenius eigenvalue $\pf$.
Let $P_M$ be its minimal polynomial:
\(
P_M = (X-\pf) \prod_{i=1}^p{(X-\alpha_i)^{m(i)}}
\).
Then, the coefficients of $A^n$ are given by:
\[
[A^n]_{ab} = c_{ab} \pf^n + \sum_{i=1}^p{P_i^{(a,b)}(n)\alpha_i^n}\, ,
\]
where the $P_i$'s are polynomials of degree $m(i)$, and $c_{ab} > 0$.
\end{lemma}

Note that the coefficients $(a,b)$ of $A^n$ gives the number of paths of length
$n$ in the diagram between some vertex $(a,k) \in \Vv_k$ and $(b,n+k)$.
This lemma states that this number is equivalent to $c_{ab} \pf^n$ when $n$
is large.

\begin{proof}
We have $c_{ab}>0$, as it is the $(a,b)$ entry of
the matrix $v_L v_R^{\tra}$ defined above.
The rest is classical, and results from the Jordan decomposition of the matrix $A$.
\end{proof}

We assume that $\Omega$ is a $d$-dimensional tiling space.
Since $\omega$ expands the distances by a factor $\lambda$, the volumes of the tiles
are dilated by $\lambda^d$. This gives the following result:
\begin{proposi}
Let $\pf$ be the Perron--Frobenius eigenvalue of $A$.
Then $\pf = \lambda^d$. In particular, $\pf > 1$.
\end{proposi}


\begin{theo}
\label{tLap09.thm-zeta}
The $\zeta$-function for the  weighted Bratteli diagram $(\partial \Bb,w)$
has abscissa of convergence $s_0=d$.
\end{theo}

\begin{proof}
We have
\begin{align*}
\zeta(s) = \sum_{\gamma \in \Pi}{w(\gamma)^s}.
\end{align*}
The quantity $w(\gamma)$ only depends on $r(\gamma)$. Furthermore,
$w(a,n)$ tends to zero like $\lambda^n = \big( \pf^{1/d} \big)^n$ when $n$
tends to infinity.
So we have:
\[
\sum_{n \in \NM)}{m^s\pf^{ns/d} \Card(\Pi_n)} \leq \zeta(s) \leq
\sum_{n \in \NM)}{M^s\pf^{ns/d} \Card(\Pi_n)},
\]
where $m$ (resp.\ $M$) is the minimum (resp.\ the maximum) of the
$w(a,1)$, $a \in \Aa$
Now, since $\Card(\Pi_n)$ grows like $\pf^n$ up to a constant (see
Lemma~\ref{tLap09.lem-rec}), we have the result.
\end{proof}

\begin{theo}
\label{tLap09.thm-measure}
The measure $\mu$ given by the Dixmier trace (Definition~\eqref{tLap09.def-Dixmier})
is well defined, and given as follows.
Let $v=(v_a)_{a \in \Aa}$ be the (right) eigenvector for $A$,
normalized such that
\( 
\sum_{e \in \Ee_0}{v_{r(e)}} = 1.
\) 
For all $\gamma \in \Pi$, let $(a,n):=r(\gamma) \in \Vv_n$.
Then:
\[
\mu ([\gamma]) = v_a \pf^{-n+1} \, .
\]
In particular, in the case of a substitution tiling, $\psi_\ast(\mu)$ is the
measure given by the \emph{frequences} of the patches, and therefore is the restriction to $\Xi$ of the unique ergodic measure on $(\Omega,\RM^d)$ (see~\cite{Sol97}).
\end{theo}

\begin{proof}
Let $\gamma = (\eps_0, e_1, \ldots, e_{n-1}) \in \Pi_n$, and $(a,n):=r(\gamma)$.
Let $\Pi_\gamma$ the subset of $\Pi$ of all paths which have $\gamma$ as a prefix.
Define:
\[
f(s) = \frac{\displaystyle\sum_{\eta \in \Pi_\gamma}{w(\eta)^s}}%
{\displaystyle\sum_{\eta \in \Pi}{w(\eta)^s}}.
\]
Then $\mu([\gamma]) = \lim_{s \downarrow d}{f(s)}$, when this limit exists.
The terms of the sum above can be grouped together: if $\eta$ is a
path of $\Pi_\gamma$, then the quantity  $w(\eta)^s$ only depends on
the length of $\eta$, say $n+k$ ($k \geq 0$), and on $r(\eta)=b \in \Aa$.
Then, if we call $N(a,b;k)$ the number of paths of length $k$ from $a$ to
$b$, we can group the sum and write:
\[
f(s) = \frac{\displaystyle\sum_{k \geq 0}{\sum_{b \in \Aa}{N(a,b;k) w(b,n+k)^s }}}%
{1+ \displaystyle\sum_{k \geq 0}{\sum_{\eps \in \Ee_0}{\sum_{b \in \Aa}{N(r(\eps),b;k) w(b,k)^s }}}}.
\]
Now, since $w(b,n)^s = \lambda^{(-n+1)s}w(b,1)^s$,
we can write:
\[
\sum_{b \in \Aa}{N(a,b;k) w(b,n+k)^s}= \lambda^{(-n-k+1)s} E_a^{\tra} A^k W(s),
\]
where $E_a$ is the vector $(\delta_a(b))_{b \in \Aa}$, and $W$ is the continuous
vector-valued function $s \mapsto (w((b,1))^s)_{b \in \Aa}$.
Similarly,
\[
\sum_{\eps \in \Ee_0}{\sum_{b \in \Aa}{N(r(\eps),b;k) w(b,k)^s}}
= \lambda^{(-k+1)s}E_{\Ee_0}^{\tra} A^k W(s),
\]
where $E_\Ee$ is the sum over $\eps \in \Ee_0$ of all $E_{r(\eps)}$.

Now, by lemma \ref{tLap09.lem-rec}, we have:
\[
E_a^{\tra} A^k W(s) = c_a (s) \pf^n + \sum_{i=1}^p{P_i(n,s)\alpha_i^n},
\]
where the $P_i$'s are polynomial in $n$ for $s$ fixed, with
$\pf > \left\vert \alpha_i \right\vert$ for all $i$, and $c_a(s)>0$.
Furthermore, the $P_i$ are continuous in $s$.
Similarly,
\[
E_{\Ee_0}^{\tra} A^k W(s) = c_\Ee (s) \pf^n + \sum_{i=1}^p{Q_i(n,s)\alpha_i^n},
\]
and since $E_{\Ee_0}$ is a linear combination of the $E_a$'s,
\begin{equation}\label{tLap09.eq-normal}
c_\Ee = \sum_{\eps \in \Ee_0}{c_{r(\eps)}}.
\end{equation}

Then, we write:
\[
\begin{split}
f(s) & = \lambda^{-ns} \frac{\displaystyle\sum_{k \geq 0}{c_a (\pf / \lambda^s)^n} + \sum_{k \geq 0}{\sum_{i=1}^p{P_i(n,s)(\alpha_i/\lambda^s)^n}} }%
{\lambda^s + \displaystyle\sum_{k \geq 0}{c_\Ee (\pf/\lambda^s)^n}+ \sum_{k \geq 0}{\sum_{i=1}^p{Q_i(n,s)(\alpha_i/\lambda^s)^n}} } \\
     & = \lambda^{-ns} \frac{\displaystyle\sum_{k \geq 0}{c_a (\pf / \lambda^s)^n} + R_1(s) }%
{\lambda^s + \displaystyle\sum_{k \geq 0}{c_\Ee (\pf/\lambda^s)^n} + R_2(s) },
\end{split}
\]
Note that the expression above is defined \emph{a priori} for $s > d$, but $R_i$ ($i=1,2$) is
defined and continuous for $s \geq d$ (the continuity results from the absolute convergence of the sum).
The remaining sums above can be computed explicitly, and we have:
\[
\begin{split}
f(s) & = \lambda^{-ns} \frac{c_a (s) + \Big(1 - (\pf / \lambda^s)^n \Big) R_1(s) }%
{c_\Ee + \Big(1 - (\pf / \lambda^s)^n \Big)( R_2(s) + \lambda^s)}.
\end{split}
\]
Then it is now clear that this expression is continuous when $s$ tends to $d$,
and has limit $\pf^{-n} c_a(d) / c_\Ee (d)$.

Let $u_a$ be defined as $c_a / c_\Ee$ for all $a$.
Let us show that $u=v$. First, show that $(u_a)_{a \in \Aa}$ (or equivalently,
$(c_a)_{a \in \Aa}$) is an eigenvector of $A$ associated with $\pf$.
We have:
\[
\begin{split}
c_a (d) & = \lim_{n \rightarrow +\infty}{\frac{E_a A^n W(d)}{\pf^n}} \\
        & = E_a L W(d),
\end{split}
\]
where $L = xy^{\tra}$, with $x$ (resp.\ $y$) an eigenvector
of $A$ (resp.\ of $A^{\tra}$) associated with $\pf$, and $\sca{x}{y} = 1$.
So $c_a$ is the $a$-coordinate of $LW(d) = \sca{y}{W(d)} x$, and so is
an eigenvector of $A$ associated to $\pf$.
Equation~\eqref{tLap09.eq-normal} now proves that $u$ has the good normalization,
and that $u=v$.
\end{proof}

\subsection{Complexity and box counting dimension}


\begin{defini}
Let $(X,d)$ be a compact metric space.
Then, the box counting dimension is defined as the following limit, when it
exists:
\[
\dim (X,d) = \lim_{t \rightarrow 0}{-\frac{\ln(N_t)}{\ln(t)}},
\]
where $N_t$ is the minimal number of balls of radius $t$ needed to
cover $X$.
\end{defini}

\begin{theo}[Thm.~2 in \cite{PB09}]
Given an ultrametric Cantor set and its associated $\zeta$-function,
let $s_0$ be the abscissa of convergence of $\zeta$.
Then:
\[
s_0 = \dim (X,d) \, .
\]
\end{theo}

This dimension can be linked to complexity for aperiodic repetitive tilings.
\begin{defini}
The complexity of a tiling $T$ is a function $p$ defined as follows:
\[
p(n) = \Card\Big\{ q \subset (T-x) \ ; \ x\in\RM^d, \  B(0,n) \subset \Supp(q)
   \text{ and }
   \forall q' \subset q, \ B(0,n) \not\subset \Supp(q') \Big\}.
\]
In other words, $p(n)$ is the number of patches $q$ of $T$ of size $n$.
\end{defini}
Note that when $T$ is repetitive (for example when $T$ is a substitution tiling)
the complexity is the same for all the tilings which are in the same tiling space.

\begin{theo}
\label{tLap09.theo-box-cplex}
Let $\Omega$ be a minimal tiling space, and $\Xi$ its canonical transversal,
endowed with the metric defined in Section~\ref{tLap09.ssect-tilings}.
Let $p$ be the associated complexity function.
Then, for this metric, the box dimension is given by:
\[
\dim(\Omega,d) = \lim_{n \rightarrow +\infty}{\frac{\ln(p(n))}{\ln(n)}}.
\]
\end{theo}

\begin{proof}
Let $N_t$ be the number of balls of diameter smaller than $t$ needed to cover
$\Xi$.
Let us first prove that for all $n \in \NM$,
\[
p(n) = N_{1/n}.
\]
Let $\Ll(n)$ be the set of all patches of size $n$, so that
$p(n) = \Card(\Ll(n))$.
Then, since all the tilings of $\Xi$ have some patch of size $n$ at the origin,
the set:
\[
\Big\{ U_q \ ; \ q \in \Ll(n) \Big\}
\]
is a cover of $\Xi$ by sets of diameter smaller than $1/n$.
So $p(n) \geq N_{1/n}$.
To prove the equality, assume we have a covering of $\Xi$ by open sets
$\{V_i \ ; \ i \in I\}$, with $\Card(I) < p(n)$, and $\diam(V_i) \leq 1/n$
for all $n$.
Then, in every $V_i$, we can find some set of the form $U_q$, $q$ a patch.
This allows us to associate some patch $q(i)$ to all $i \in I$.
We claim that for all such $q(i)$, $B(0,n) \subset \Ext (q(i))$, where $\Ext(q)$
is the maximal unambiguous extension of $p$, as defined in \ref{tLap09.def-unambiguous}.
Indeed, if it were not the case, $\diam(U_q)$ would be smaller than $n$.
Therefore, by restriction, to each $i \in I$, we can associate a patch $q'(i)$
of size $n$, such that $U_{q'(i)}  \subset V_i$.
Since the $\{V_i\}_{i \in I}$ cover $\Xi$, all patches of size $n$ are obtained
this way, and $p(n) \leq \Card(I)$.
Since this holds for all cover, $p(n) \leq N_{1/n}$.

Now, $N_t$ is of course an increasing function of $t$.
Therefore,
\[
N_{t-1} \leq p\big([1/t]\big) \leq N_t,
\]
and so:
\[
-\frac{N_{t-1}}{\ln(t)} \leq \frac{p\big([1/t]\big)}{\ln([1/t])} \leq -\frac{N_{t}}{\ln(t-1)}.
\]
Letting $t$ tend to zero proves the theorem.
\end{proof}

\begin{coro}
\label{tLap09.coro-cplex-subs}
Let $\Xi$ be the transversal of a minimal aperiodic tiling space, with a
complexity function which satisfies \(C_1 n^\alpha \leq p(n) \leq C_2 n^\alpha\)
for $C_1, C_2, \alpha > 0$.
Then:
\[
\dim(\Xi,d) = \alpha.
\]
\end{coro}

Let now consider a tiling space $\Omega$ associated to a substitution $\omega$.
Let $\Bb$ be the weighted Bratteli diagram associated with it.
We proved in Section~\ref{tLap09.ssect-zetameasure} that for a Bratteli diagram associated with a
substitution tiling of dimension $d$, the abscissa of convergence is
exactly $d$.
It is furthermore true that the box dimension of the transversal $\Xi$
with the usual metric is $d$; it results from the invariance of the box dimension under
bi-Lipschitz equivalence, which is proved in \cite[Ch.~2.1]{Fal90}.

Therefore, we can deduce the following result:
\begin{coro}
\label{tLap09.thm-complex}
Let $\Omega$ be a substitution tiling space satisfying our conditions.
Then there exists a function $\nu$ such that:
\[
p(n) = n^{\nu(n)}, \quad 
  \textrm{with } \lim_{n \rightarrow +\infty}{\nu(n)} = d \, .
\]
Equivalently, for all $\epsilon > 0$, there exists $C_1,C_2>0$ such that for all
$n$ large enough,
\[
C_1 n^{d-\epsilon} \leq p(n) \leq C_2 n^{d+\epsilon} \,.
\]
\end{coro}

This result is actually weaker than what we can actually expect:
in fact, there exist $C_1,C_2 > 0$ such that
\[
C_1 n^d \leq p(n) \leq C_2 n^d \, .
\]
The upper bound was proved by Hansen and Robinson for self-affine tilings
(see~\cite{Rob04}).
The lower bound can be proved by direct analysis for substitution tilings.
It would also result from the conjecture that any $d$-dimensional tiling with
low complexity (which means $p(n)/n^d$ tends to zero) has at least one period
(see~\cite{LP03}).

However, it is still interesting to see how the apparently abstract fact that
the abscissa of convergence $s_0$ equals the dimension gives in fact a result
on complexity.

\section{Laplace--Beltrami operator}
\label{tLap09.sect-LBop}

Let $\Bb$ be a weighted Bratteli diagram.
The Dixmier trace \eqref{tLap09.eq-Dixmier} induces a probability measure $\nu$ on the set $E$ of choice functions (see \cite{PB09} section 7.2.).
The following is \cite{PB09} Theorem 4.

\begin{proposi} 
\label{tLap09.prop-Dirichlet}
For all $s \in \RM$ the bilinear form on $L^2(\partial \Bb, d\mu)$ given by
\begin{equation}
\label{tLap09.eq-Dirichlet}
Q_s(f,g) = \frac{1}{2} \int_E \TR
\bigl(
|D|^{-s} [D,\pi_\tau(f)]^\ast [D, \pi_\tau(g)]
\bigr) \; d\nu(\tau)\,,
\end{equation}
with dense domain 
\(\text{\rm Dom} \, Q_s = 
\langle 
\chi_\gamma \, : \, \gamma \in \Pi
\rangle\), is a closable {\em Dirichlet form}.
\end{proposi}

The classical theory of Dirichlet forms \cite{Fuk80} allows to identify $Q_s(f,g)$ with \(\langle f, \Delta_s g \rangle\) for a non-positive definite self-adjoint operator $\Delta_s$ on $L^2(\partial \Bb, d\mu)$ which is the generator of a Markov semi-group.
We have \( \text{\rm Dom}\, Q_s \subset \text{\rm Dom}\,\Delta_s \subset \text{\rm Dom}\,\tilde{Q}_s\) where $\tilde{Q}_s$ is the smallest closed extension of $Q_s$. 
The following is taken from \cite{PB09} section 8.3.

\begin{theo}
\label{tLap09.thm-PBDelta}
The operator $\Delta_s$ is self-adjoint and has pure point spectrum.
\end{theo}

Following Pearson and Bellissard we can calculate $\Delta_s$ explicitly on characteristic functions of cylinders. For a path \( \eta \in \Pi\) let us denote by $\ext(\eta)$ the set of ordered pairs of distinct edges $(e,e')$ which can extend $\eta$ one generation further.

\begin{subequations}
\label{tLap09.eq-LapExplicit}
\begin{align}
\Delta_s \chi_\gamma = 
- \sum_{k=0}^{|\gamma|-1} \frac{1}{G_s(\gamma_k)}
\Bigl(
( \mu[\gamma_{k}]-\mu[\gamma_{k+1}] ) \chi_\gamma 
- \mu[\gamma] ( \chi_{\gamma_k}-\chi_{\gamma_{k+1}} )
\Bigr)
\label{tLap09.eq-Lap}\\
\text{with} \quad G_s(\eta)= \frac{1}{2} \ \text{\rm diam}[\eta]^{2-s} 
\sum_{(e,e')\in \ext(\eta)} \mu[\eta\cdot e]\mu[\eta\cdot e'] \quad \quad
\label{tLap09.eq-G}
\end{align}
\end{subequations}

Note that the term \(\chi_{\gamma_k}-\chi_{\gamma_{k+1}}\) is the characteristic function of all the paths which coincide with  $\gamma$ up to generation $k$ and differ afterwards, {\it i.e.} all paths which split from $\gamma$ at generation $k$.
And \(\mu[\gamma_{k}]-\mu[\gamma_{k+1}]\) is the measure of this set.
 
We now state the main theorem which gives explicitly the full spectrum of $\Delta_s$.

\begin{theo}
\label{tLap09.thm-SpectrumLap}
The spectrum of $\Delta_s$ is given by the following.
\begin{enumerate}[(i)]

\item $0$ is a single eigenvalue with eigenspace 
\(\langle 1 = \chi_{\partial \Bb} \rangle\).

\item $\lambda_0 = \frac{1}{G_s(\circ)}$ is eigenvalue with eigenspace 
\(E_0 = \langle \frac{1}{\mu[\eps]} \chi_\eps - \frac{1}{\mu[\eps']} \chi_{\eps'} \, : \, \eps,\eps' \in \Ee_0\,, \eps\ne \eps'\rangle\) of dimension \(\dim E_0 = n_0 -1\), where $n_0$ is the cardinality of $\Ee_0$.

\item For $\gamma \in \Pi$, 
\begin{equation}
\label{tLap09.eq-eigenfn}
\lambda_\gamma = \sum_{k=0}^{|\gamma|-1} 
\frac{\mu[\gamma_{k+1}]-\mu[\gamma_k]}{G_s(\gamma_k)} - \frac{\mu[\gamma]}{G_s(\gamma)}
\end{equation}
is eigenvalue with eigenspace
\begin{equation}
\label{tLap09.eq-eigenspace}
E_\gamma = 
\Bigl< 
\frac{1}{\mu[\gamma\cdot e]} \chi_{\gamma\cdot e} - 
\frac{1}{\mu[\gamma\cdot e']} \chi_{\gamma\cdot e'} \, : \,
(e,e') \in \ext(\gamma)
\Bigr>
\end{equation}
of dimension \(\dim E_\gamma = n_\gamma -1\), where $n_\gamma$ is the number of edges extending $\gamma$ one generation further.

\end{enumerate}
\end{theo}
\begin{proof}
The formula for the eigenvalues is calculated easily noticing that $\Delta_s \chi_{\gamma\cdot e}$ and $\Delta_s \chi_{\gamma \cdot e'}$ only differ by the last term in the sum in equation \eqref{tLap09.eq-Lap}.

The spectrum of $\Delta_s$ is always the closure of its set of eigenvalues (whatever its domain may be).
Hence we do not miss any of it by restricting to characteristic functions.

We now show that all the eigenvalues of $\Delta_s$ are exactly given by the $\lambda_\gamma$.
It suffices to check that the restriction of $\Delta_s$ to $\Pi_n$ has exactly $\dim \Pi_n$ eigenvalues (counting multiplicity).
Notice that $\chi_\gamma$ is the sum of $\chi_{\gamma \cdot e}$ over all edges $e$ extending $\gamma$ one generation further.
Hence an eigenfunction in $E_\gamma$ for $\gamma \in \Pi_k$ can be written as a linear combination of characteristic functions of paths in $\Pi_m$, for any $m > k$.
The number of eigenvalues $\lambda_\gamma$ for $\gamma \in \Pi_{n-1}$ is \( \sum_{\gamma \in \Pi_{n-1}} \dim E_\gamma = \sum_{\gamma \in \Pi_{n-1}} (n_\gamma -1) = \dim \Pi_n - \dim \Pi_{n-1}\).
So the number of eigenvalues $\lambda_\gamma$ for $\gamma \in \Pi_k$ for all $1\le k \le n-1$ is \( \sum_{k=1}^{n-1}\bigl( \dim \Pi_{k+1} - \dim \Pi_k \bigr) = \dim \Pi_n - \dim \Pi_1\).
And counting $0$ and $\lambda_0$ adds up \((\dim \Pi_1 - 1) + 1\) to make the count match.
\end{proof}

\begin{rem}
\label{tLap09.rem-spectrumPB}
{\em As noted in Remark \ref{tLap09.rem-tree}, our formalism with Bratteli diagrams includes as a special case the approach of Pearson--Bellissard for weighted Cantorian Michon trees. Hence Theorem \ref{tLap09.thm-SpectrumLap} gives also the spectrum and eigenvectors of their Laplace--Beltrami operators.
}
\end{rem}

The eigenvectors \eqref{tLap09.eq-eigenspace} are very simple to picture.
Given a path $\gamma$, and two extensions \((a,b) \in \ext(\gamma)\), an eigenvector is the difference of their characteristic functions weighted by their measures.
See Figure \ref{tLap09.fig-EfnFibo} in section \ref{tLap09.sect-intro} for an example for the Fibonacci diagram.

\section{Cuntz--Krieger algebras and applications}
\label{tLap09.sect-CK}

We now consider stationary Bratteli diagrams. 
We use the self-similar structure to further characterize the operator $\Delta_s$ and its spectrum.

\subsection{Cuntz--Krieger algebras}
\label{tLap09.ssect-CKalg}

Let $\Bb$ be a stationary Bratteli diagram.
Let $A$ be its Abelianization matrix.
Let us denote by $\Ee_0$ the set of edges $\eps$ linking the root to generation $1$, and
by $\Emod$ the set of edges linking two generations (excluding the root).
Let \(\tilde{A} = \bigl( \tilde{a}_{ef}\bigr)_{ e,f \in \Emod}\) be the square matrix with entries $\tilde{a}_{ef}=1$ if $e$ can be composed with (or followed down the diagram by)  $f$ and $\tilde{a}_{ef}=0$ else.
There is an associated Bratteli diagram $\tilde{\Bb}$ with Abelianization matrix $\tilde{A}$, which is ``dual'' to $\Bb$ in the sens that its vertices corresponds to the edges of $\Bb$ and its edges to the adjacencies of edges in $\Bb$.
Note that, because the entries of $\tilde{A}$ are zeros or ones, all the edges in $\tilde{\Bb}$ are simple.

\begin{rem}
\label{tLap09.rem-vertices}
{\em
If $\Bb$ has only simple edges, we can simply take $\tilde{\Bb} = \Bb$ and $\tilde{A}=A$.
}
\end{rem}

The Cuntz-Krieger algebra $\Oo_{\tilde{A}}$, is the \Cs generated by the partial isometries $U_e, e \in \Emod$ (on a separable, complex, and infinite dimensional Hilbert space $\Hh$) that satisfy the following relations.

\begin{equation}
\label{tLap09.eq-CKalg}
\Oo_{\tilde{A}} = C^\ast
\bigl< \, U_e, U_e^\ast \,, e\in \Emod \ \ \vert \ U_e U_e^\ast, U_e^\ast U_e \in \Pp(\Hh)\,, \ 
U_e^\ast U_e = \sum_{f\in\Emod} \tilde{A}_{ef} \, U_f U_f^\ast \, 
\bigr>\,,
\end{equation}
where $\Pp(\Hh)$ denotes the set of projections in $\Hh$: \(p \in \Pp(\Hh) \iff p^2=p^\ast=p\).

\begin{figure}[!h]
\[
\xymatrix{
& \circ \ar@{.}[dl] \ar@{-}[dr]^{\eps} &                         
&& & \circ \ar@{-}[dl]_{\eps'} \ar@{.}[dr] & 
&& & \circ \ar@{-}[dl]_{\eps'} \ar@{.}[dr] & \\
\bullet \ar@{:}[dd]^{e_2}_{e_1} \ar@{.}[ddrr]^(.3){e_3} & & \bullet \ar@{- }[ddll]_(.3){e_4} \ar@{.}[dd]^{e_5}
&&  \bullet\ar@{:}[dd] \ar@{-}[ddrr] & & \bullet \ar@{.}[ddll]  \ar@{.}[dd]
&&  \bullet\ar@{:}[dd] \ar@{-}[ddrr]^(.3){U_{e_4}^\ast \gamma} & & \bullet\ar@{.}[ddll]  \ar@{.}[dd] \\
& & && & & && & & \\
\bullet\ar@{:}[dd] \ar@{-}[ddrr]_(.3){\gamma} & & \bullet\ar@{.}[ddll] \ar@{.}[dd]              
&&  \bullet\ar@{:}[dd] \ar@{.}[ddrr] & & \bullet\ar@{-}[ddll]  \ar@{.}[dd]
&& \bullet\ar@{:}[dd] \ar@{.}[ddrr] & & \bullet\ar@{.}[ddll] \ar@{.}[dd] \\
& & && & & && & & \\
\bullet\ar@{:}[dd] \ar@{.}[ddrr] & & \bullet\ar@{.}[ddll] \ar@{.}[dd]       
&&  \bullet\ar@{:}[dd] \ar@{-}[ddrr]^(.3){U_{e_3} \gamma} & & \bullet\ar@{.}[ddll]  \ar@{.}[dd] 
&& \bullet\ar@{:}[dd] \ar@{.}[ddrr] & & \bullet\ar@{.}[ddll] \ar@{.}[dd] \\
& & && & & && & & \\
\bullet& &\bullet && \bullet& &\bullet && \bullet& &\bullet 
}
\]
\caption[haha]{{\small Example of some Cuntz--Krieger operators acting on finite paths.}}
\label{tLap09.fig-CK}
\end{figure}
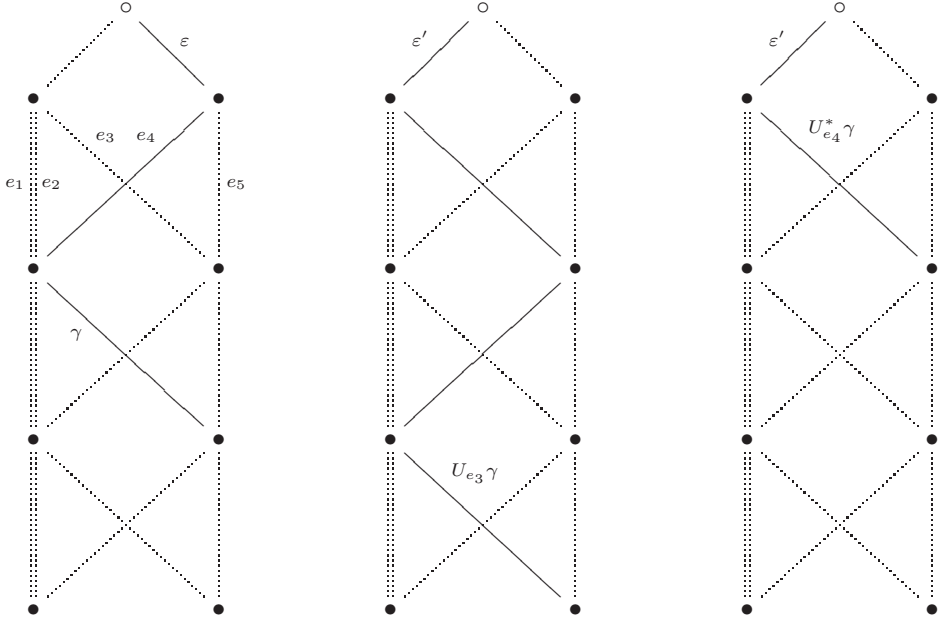

By abuse of notation we write the basis elements of $l^2(\Pi \setminus \Pi_1)$ as $\gamma \in \Pi \setminus \Pi_1$.
The Cuntz-Krieger algebra $\Oo_{\tilde{A}}$ is represented on $l^2(\Pi \setminus \Pi_1)$ as follows.
\begin{subequations}
\label{tLap09.eq-actionCK-l2}
\begin{align}
U_e (\varepsilon, e_{1}, e_{2}, \cdots) = 
\left\{ 
\begin{array}{cl}
(\varepsilon', e, e_{1}, e_{2}, \cdots ) & \text{\rm if} \ A_{e e_1} = 1 \\
0 & \text{\rm else}
\end{array}
\right.
\label{tLap09.eq-Ui}\\
U_e^\ast (\varepsilon, e_{1}, e_{2}, \cdots) = 
\left\{ 
\begin{array}{cl}
(\varepsilon', e_{2}, e_{3}, \cdots ) & \text{\rm if} \ e_{1} = e \\
0 & \text{\rm else}
\end{array}
\right.
\quad \
\label{tLap09.eq-Ui*}
\end{align}
\end{subequations}
where $\varepsilon' \in \Ee_0$ stands for the (possibly) new edges linking to the root, as illustrated in the case of the Penrose substitution below.
The orientation of $\varepsilon'$ however is taken to be the same as that of $\varepsilon$:
if \( \eps=(s(e_1),g)\) then we have \(\eps'=(s(e),g)\) (see Remark \ref{tLap09.rem-group}).
In other words we require that
\(\psi  (\eps', e, e_{1}, \cdots)  = \omega \circ \psi  (\eps, e_{1}, \cdots) + x\) for some $x \in \RM^d$, and where \(\psi: \partial \Bb \rightarrow \Xi\) is the homeomorphism of Proposition \ref{tLap09.prop-homeo}. 
See Figure \ref{tLap09.fig-CK} for some examples.

There is an induced ``action'' on \( \text{\rm Dom} \Delta_s \subset L^2(\partial \Bb, d\mu)\) defined as follows:
\begin{equation}
\label{tLap09.eq-actionCK-L2}
U_e^{(\ast)} \chi_\gamma = 
\left\{
\begin{array}{cl}
\chi_{U_e^{(\ast)} \gamma} & \text{\rm if} \ U_e^{(\ast)} \gamma \neq 0\\
0 & \text{\rm else}
\end{array}
\right.
\end{equation}
for paths in $\Pi_{n\ge 2}$ and by linearity for shorter paths, using the relations \(\chi_\gamma = \sum_e \chi_{\gamma\cdot e}\) (the sum running over all edges $e$ extending the path one generation further).

If $\varphi \in E_\gamma$, for $\gamma \in \Pi_{n\ge 2}$, is an eigenfunction of $\Delta_s$, then we see from Equation~\eqref{tLap09.eq-actionCK-L2} that $U_e \varphi$, if not zero, is another eigenfunction.
Since the diagram $\Bb$ is stationary, we have the following obvious fact.
\begin{equation}
\label{tLap09.eq-actionCK-efn}
 U_e \ E_\gamma = 
\left\{
\begin{array}{cl}
E_{U_e \gamma} & \text{\rm if} \ U_e \gamma \neq 0 \\
\{ 0 \} & \text{\rm else}\,.
\end{array}
\right.
\end{equation}
This allows to define an induced ``action'' of $\Oo_{\tilde{\Aa}}$ on the set of eigenvalues of $\Delta_s$ as follows.
\begin{equation}
\label{tLap09.eq-actionCK-eval}
u_e(\lambda_\gamma) = 
\left\{
\begin{array}{cl}
\lambda_{U_e \gamma} & \text{\rm if} \ U_e \gamma \neq 0 \\
 0  & \text{\rm else}\,.
\end{array}
\right.
\end{equation}
Those maps $u_e$ are calculated as shown below.
\begin{lemma}
\label{tLap09.lem-CK-eval}
Let $\pf$ be the Perron--Frobenius eigenvalue of $A$.
Let \(\gamma=(\varepsilon,e_1, e_2, \cdots ) \in \Pi_{n\ge 2}\).
Then, if $U_e \gamma \neq 0$, we have
\begin{equation}
\label{tLap09.eq-CK-eval-expl}
u_e (\lambda_\gamma) =
\pf^{(d+2-s)/d} 
\Bigl( 
\lambda_\gamma - \frac{ \mu[\eps]-\mu[\circ] }{ G_s(\circ) } 
\Bigr)
+ \frac{ \mu[\eps']-\mu[\circ] }{ G_s(\circ) }
+ \frac{ \mu[\eps' e] -\mu[\eps'] }{ G_s(\eps') }
\end{equation}
\end{lemma}
\begin{proof}
By equation \eqref{tLap09.eq-eigenfn} we have
\begin{equation}
\label{tLap09.eq-lem-CK-eval-1}
\lambda_{U_e\gamma} = - \sum_{k=0}^{|U_e\gamma|-1} 
 \frac{\mu[(U_e\gamma)_{k}]-\mu[(U_e\gamma)_{k+1}]}{G_s((U_e\gamma)_k)} - \frac{\mu[U_e\gamma]}{G_s(U_e\gamma)} \,.
\end{equation}
The terms for $k=0,1,$ in the above sum give the last two terms in equation \eqref{tLap09.eq-CK-eval-expl}.
For all $k\ge 2$, \(\mu[(U_e\gamma)_k=\pf^{-1}\mu[\gamma_{k-1}]\), and \(G_s((U_e\gamma)_k) = \pf^{(s-d)/d-2}G_s(\gamma_{k-1})\) (see Theorem \ref{tLap09.thm-measure} for the rescalling of the measures).
Hence the rest of the sum, over \( k=2, \cdots |U_e\gamma|-1\) (and the last term) in equation \eqref{tLap09.eq-lem-CK-eval-1}, rescale by a factor $\pf^{(d-s+2)/d}$ to the sum over \( k=1, \cdots |\gamma|-1\) (and the last term) for the eigenvalue $\lambda_\gamma$ (equation \eqref{tLap09.eq-eigenfn}).
We then add the contribution of the root, {\it i.e.} the term for $k=0$, to get equation \eqref{tLap09.eq-CK-eval-expl}.
\end{proof}

Note that the maps $u_e$ are affine, with constant terms, written $\beta_e$, that only depend on $e$.
We will write for now on
\begin{equation}
\label{tLap09.eq-ui(eval)}
u_e(\lambda_\gamma) = \Lambda_s \lambda_\gamma + \beta_e \,, 
\quad \quad \text{\rm with} \quad
\Lambda_s = \pf^{(d+2-s)/d}\,.
\end{equation}

The eigenelements of $\Delta_s$ corresponding to $E_0$, and $E_{\eps}, \eps \in \Ee_0$, are immediate to calculate explicitly from equations \eqref{tLap09.eq-eigenfn} and \eqref{tLap09.eq-eigenspace}.
We can therefore calculate explicitly all other eigenelements by action of $\Oo_{\tilde{\Aa}}$ on those corresponding to the $E_{\eps}, \eps \in \Ee_0$.
We summarize this in the following.

\begin{proposi}
\label{tLap09.prop-ui(eval)-explicit}
For \(\gamma = (\varepsilon,e_{1}, e_{2}, \cdots e_{n}) \in \Pi_{n\ge 1}\), set
\(U_\gamma = U_{e_1} U_{e_2} \cdots U_{e_n}\) and \( u_\gamma = u_{e_1} \circ u_{e_2} \circ \cdots u_{e_n}\). 
For any \(\gamma \in \Pi_{n\ge 1}\), we have \(E_\gamma = U_\gamma E_{\varepsilon'}\), and
\begin{equation}
\label{tLap09.eq-CK-spectrum-eval}
 \lambda_\gamma = u_\gamma (\lambda_{\varepsilon'}) = 
\, \Lambda_s^n \lambda_{\varepsilon'} \ + \  \sum_{j=1}^n \, \Lambda_s^{j-1} \,\beta_{e_j} \,.
\end{equation}
\end{proposi}

\subsection{Bounded case}
\label{tLap09.ssect-bdcase}

We consider here the case \(s>d+2\).
We show that $\Delta_s$ is bounded and characterize the boundary of its spectrum.

\begin{proposi} 
\label{tLap09.prop-bdcase}
For  \(s>d+2\), \( \Delta_s\) is a {\em bounded}, and we have
\[
\| \Delta_s \|_{\Bb (  L^2 (\partial \Bb, d\mu) )} \le \ c
\ \frac{1}{1- \pf^{(d+2-s)/d}} \,,
\]
with \( c = \max \bigl\{ \max_{\eps \in \Ee_0} |\lambda_\eps|, \max_{e\in\Emod}|\beta_e| \bigr\}\).
\end{proposi}
\begin{proof}
By Proposition \ref{tLap09.prop-ui(eval)-explicit}, equation \eqref{tLap09.eq-CK-spectrum-eval}, we see that for $\gamma \in \Pi_n$ we have.
\[
|\lambda_\gamma| \le 
\max \bigl\{ \max_{\eps \in \Ee_0} |\lambda_\eps|, \max_{e\in\Emod}|\beta_e|
\bigr\}
\ \sum_{j=0}^n \Lambda_s^j \,.
\]
From equation \eqref{tLap09.eq-ui(eval)} we have \(\Lambda_s = \pf^{(d-s+2)/d} < 1 \),
therefore the above geometric sum converges, and is bounded for all $n$ by its sum \( 1/(1-\Lambda_s)\).
\end{proof}

We define the $\omega$-spectrum of $\Delta_s$ as
\[
\Sp_\omega (\Delta_s) = \bigcap_{n \in \NM} \overline{\Sp(\Delta_s) \setminus \Sp \bigl( \Delta_s \vert_{\Pi_n} \bigr)} \,.
\]
In our case here, this is the boundary of the (pure point) spectrum of $\Delta_s$.
Under some conditions on $A$ and the $\beta_e, e \in \Emod$, and for \(s>d+2\) large enough, one can show that \(\Sp_\omega(\Delta_s)\) is homeomorphic to $\Xi$, and that this homeomorphism is {\em  H\"older} \cite{JS09}.

\subsection{Weyl asymptotics}
\label{tLap09.ssect-Weyl}

The following theorem justifies calling $\Delta_s$ a Laplace--Beltrami operator.
Indeed, for $s=s_0=d$, Theorem \ref{tLap09.thm-Weyl} shows that the number of eigenvalues of $\Delta_d$ of modulus less that $\lambda$ behaves like $\lambda^{d/2}$ when \( \lambda \rightarrow \infty\), which is the classical Weyl asymptotics for the Laplacian on a compact $d$-manifold.

\begin{theo}
\label{tLap09.thm-Weyl}
Let 
\( \Nn_s(\lambda) = \text{\rm Card}
\bigl\{
\lambda' \;\text{eigenvalue of} \; \Delta_s \ : \
|\lambda'| \le \lambda
\bigr\}
\). 
For $s<d+2$, we have the following {\rm Weyl asymptotics}
\begin{equation}
\label{tLap09.eq-Weyl}
 c_- \lambda^{d/(d-s+2)} \le \Nn_{s}(\lambda) \le  c_+ \lambda^{d/(d-s+2)}\,,
\end{equation}
as \(\lambda \rightarrow + \infty\), for some constants $0<c_-<c_+$.
\end{theo}
\begin{proof}
By Proposition \ref{tLap09.prop-ui(eval)-explicit}, equation \eqref{tLap09.eq-CK-spectrum-eval}, there exists constants $x_\pm, y_\pm>0$, such that for all $\gamma \in \Pi_n$ we have
\( x_- \Lambda_s^n + y_- \le |\lambda_\gamma| \le  x_+ \Lambda_s^n + y_+\).
For $s<d+2$, $\Lambda_s>1$ (see equation \eqref{tLap09.eq-ui(eval)}), so there is an integer $k>0$ (independent of $n$) such that
\( x_+ \Lambda_s^n + y_+ \le  x_- \Lambda_s^{n+k} + y_-\).
Hence for all $\gamma \in \Pi_l, \, l\le n$, we have
\( |\lambda_\gamma| \le  x_+ \Lambda_s^n + y_+\), 
and for all $\gamma \in \Pi_l, \, l\ge n+k$, we have 
\( |\lambda_\gamma| \ge  x_- \Lambda_s^{n+k} + y_-\).
Therefore we get the inequalities: 
\( \Card \, \Pi_n \le  \Nn(x_+ \Lambda_s^n + y_+) \le \Card \, \Pi_{n+k}\).
There are constants $c_1>c_2>0$ such that for all $l\in \NM$, 
\(c_1 \pf^l \le \Card \, \Pi_l \le c_2 \pf^n\), so that we get
\( c_1 \pf^n \le  \Nn(x_+ \Lambda_s^n + y_+) \le c_2 \pf^{n+k} = c_3 \pf^n\).
We substitute $\Lambda_s$ from equation \eqref{tLap09.eq-ui(eval)} to complete the proof.
\end{proof}

\subsection{Seeley equivalent}
\label{tLap09.ssect-Seeley}
For the case $s=s_0=d$ we give an equivalent to the trace \(\TR \bigl( e^{t\Delta_d} \bigr)\), for $s=s_0=d$, as $t\downarrow 0$.
The behavior of \(\TR \bigl( e^{t\Delta_d} \bigr)\) like $t^{-d/2}$ as $t\downarrow 0$ is in accordance with the leading term of the classical Seeley expansion for the heat kernel on a compact $d$--manifold.

\begin{theo}
\label{tLap09.thm-Seeley}
There exists constants $c_+\ge c_- >0$, such that as $t \downarrow 0$
\begin{equation}
\label{tLap09.eq-Seeley}
c_- t^{-d/2} \le \TR \bigl( e^{t\Delta_d} \bigr) \le c_+ t^{-d/2} \,.
\end{equation}
\end{theo}
\begin{proof}
Let $P_\gamma$ be the spectral projection (onto $E_\gamma$) for $\gamma \in \Pi_{n\ge 1}$, and $P_0$ that on $E_0$.
The trace reads
\begin{equation}
\label{tLap09.eq-trace1}
\TR \bigl( e^{t\Delta_s} \bigr) = 1 + e^{\lambda_0 t} \TR(P_0) + 
\sum_{n=1}^\infty \sum_{\gamma \in \Pi_n} e^{t \lambda_\gamma} \TR(P_\gamma) \,.
\end{equation}
Now $\TR(P_\gamma) = n_\gamma -1$, with $n_\gamma$ the number of possible extensions of $\gamma$ one generation further (see equation \eqref{tLap09.eq-eigenspace} in theorem \ref{tLap09.thm-SpectrumLap}), and $\TR(P_0) = n_0 -1$.
Since the Bratteli diagram of the substitution is stationary, the integers $n_\gamma$ are bounded, so there are $p_-, p_+>0$, such that for all $\gamma \in \Pi$ we have:
\begin{equation}
\label{tLap09.eq-bdP}
p_- \le \TR(P_\gamma) \le p_+ \,.
\end{equation}
By equation \eqref{tLap09.eq-CK-spectrum-eval}, the $\lambda_\varepsilon, \beta_e$, being bounded, there exists $\lambda_-, \lambda_+>0$, such that for all $\gamma \in \Pi$ we have: 
\begin{equation}
\label{tLap09.eq-bdlambda}
\lambda_- \pf^{nd/2} \le |\lambda_\gamma| \le \lambda_+ \pf^{nd/2}\,.
\end{equation}
The cardinality of $\Pi_n$ grows like $\pf^n$ so there are $\pi_-,\pi_+>0$, such that for all $n\ge 0$ we have:
\begin{equation}
\label{tLap09.eq-bdPi}
\pi_- \pf^n \le | \Pi_n | \le \pi_+ \pf^n\,.
\end{equation}
We substitute inequalities \eqref{tLap09.eq-bdP},  \eqref{tLap09.eq-bdlambda}, and \eqref{tLap09.eq-bdPi} into equation \eqref{tLap09.eq-trace1} to get
\begin{equation}
\label{tLap09.eq-trace2}
1+ p_- \pi_- \sum_{n=0}^\infty \pf^n e^{-t \lambda_+ \pf^{nd/2}}
\le \TR \bigl( e^{t\Delta_s} \bigr) \le
1+ p_+ \pi_+ \sum_{n=0}^\infty \pf^n e^{-t \lambda_- \pf^{nd/2}} \,.
\end{equation}
Set \(N_t = d \log(1/t)/(2 \log(\pf))\), and split the above sums into two parts: the sum over $n  < N_t$, and  the remainder.
For the finite sum we have:
\begin{equation}
\label{tLap09.eq-sum1.1}
\sum_{n=0}^{N_t -1} \pf^n e^{-t \lambda_\pm \pf^{nd/2}}  =
t^{-d/2} \Lambda^{-N_t} \sum_{n=0}^{N_t-1} \pf^n   e^{-\lambda_\pm \pf^{(n-N_t)d/2}} \,, 
\end{equation}
where we have used $\pf^{N_td/2} = 1/t$.
With the inequalities \( e^{-\lambda_\pm} \le e^{-\lambda_\pm \pf^{(n-N_t)d/2}} \le 1\), the above sum on the right hand side of \eqref{tLap09.eq-sum1.1} is bounded by the geometric series \( \sum_{n=0}^{N_t-1} \pf^n = (\pf^{N_t}-1)/(\pf-1)\).
Multiplying by  $t^{-d/2} \Lambda^{-N_t}$ we get the inequalities:
\begin{equation}
\label{tLap09.eq-sum1.2}
t^{-d/2}  c'_- 
\le \sum_{n=0}^{N_t -1} \pf^n e^{-t \lambda_\pm \pf^{nd/2}} \le
t^{-d/2} c'_+ \,,
\end{equation}
for some constants $c'_-, c'_+>0$, and $t$ small enough.

For the remainder of the sums in \eqref{tLap09.eq-trace2} we have
\begin{equation}
\label{tLap09.eq-sum2.1}
\sum_{n=N_t}^\infty \pf^n e^{-t \lambda_\pm \pf^{nd/2}} = 
\pf^{N_t} \sum_{m=0}^\infty \pf^m   e^{- \lambda_\pm \pf^{md/2}} =
t^{-d/2} c''_\pm \,,
\end{equation}
where $c''_\pm>0$ is the sum of the absolutely convergent series.
We put together inequalities \eqref{tLap09.eq-sum1.2} and equation \eqref{tLap09.eq-sum2.1} into inequalities \eqref{tLap09.eq-trace2} to complete the proof.
\end{proof}

\subsection{Eigenvalues distribution}
\label{tLap09.ssect-Evdist}

We now restrict to the case $s=s_0=d$, so that $\Lambda_d = \pf^{2/d}$ in equation \eqref{tLap09.eq-ui(eval)}.
We also suppose here that the weights are simply given by the $d$-th roots of the measures: \( \diam [\gamma] = \mu[\gamma]^{1/d}\).

Let us consider the cases $d=1$ and $d=2$ first. 
One sees from equations \eqref{tLap09.eq-G}, \eqref{tLap09.eq-eigenfn}, and \eqref{tLap09.eq-CK-eval-expl}, that the \(\lambda_\eps, \eps\in \Ee_0\), and \(\beta_e, e\in \Emod\), are rational functions of $\pf$. 
Hence by equation \eqref{tLap09.eq-CK-spectrum-eval} all the eigenvalues of $\Delta_d$ belong to the field \(\QM[\pf]\).
The maps $u_e$ in equation \eqref{tLap09.eq-ui(eval)} become affine maps in this field, whose linear part is the multiplication by $\pf^{2/d}$ ($=\pf^2$ or $\pf$).

Let \( P_A(X) = X^r + a_{r-1}X^{r-1} + \cdots a_1 X + a_0\) be the minimal polynomial of $A$.
The field $\QM[\pf]$ is isomorphic to $\QM^r$, so that the eigenvalues of $\Delta_d$ can be represented by points in  $\QM^r$.
The $\lambda_\eps, \beta_e$, are identified with fixed vectors $\vec{\lambda}_\eps, \vec{\beta}_e \in\QM^r$, and the maps $u_e$ become affine maps in $\QM^r$.
The multiplication by $\pf^{2/d}$ (linear parts of the maps $u_e$) is implemented in the basis \(\QM \oplus \pf^{2/d} \QM \oplus \cdots  \pf^{2(r-1)/d} \QM\) by the companion matrix of $P_A$:
\[
\left[
\begin{array}{ccccrc}
0 &        &&   &    & - a_0   \\
1 & \ddots &&   &    & \vdots  \\
  & \ddots && 0 &    & -a_{r-2} \\
  &        && 1 &    & - a_{r-1}\\
\end{array}
\right] \,.
\]
Note that it is equivalent to $A$ and therefore has the same Perron--Frobenius eigenvector.
Since there are finitely many $\vec{\beta}_e, \vec{\lambda}_\eps \in \QM^r$, upon multiplication by a large enough integer, we can actually represent the eigenvalues by points of the lattice $\ZM^r$.
It is then natural to ask which points correspond to eigenvalues of $\Delta_s$ and if one can characterize this set.
The answer is surprisingly simple:

\noindent {\em Eigenvalues correspond to points of integer coordinates in $\QM^r$ that stay within a bounded strip to the Perron--Frobenius eigenline of $A$.}

See section \ref{tLap09.ssect-fibo} for an example, and Figure \ref{tLap09.fig-fibo} in section \ref{tLap09.sect-intro}.

For $d\neq 1,2$ a similar result hold, but the points can no longer be chosen to have integer coordinates.
In order to prove this in general -- for $d$ generic --  we consider a quotient ring of $\QM[\pf]$ where we can implement the multiplication by $\pf^{2/d}$.
Define $Q_A(x) = P_A(x^{d/2})$ if $d$ is even, and let $C$ be its companion matrix.
If $d$ is odd define $Q_A(x) = P_A(x^d)$, and let $C$ be the square of its companion matrix.
Let us denote by $d'$ either $d/2$ if $d$ is even, or $d$ if it is odd.
The matrix $C$ implements the multiplication by $\pf^{d'}$ on the ring $\QM[X] / Q_A$, which we embed in the vector space \( V = \RM \oplus \pf \RM \oplus \pf^2 \RM \oplus \cdots  \pf^{rd'-1} \RM\).

Because of the factor \( \text{\rm diam}[\eta]^{2-s} = \mu[\eta]^{(2-d)/d}\) in equation \eqref{tLap09.eq-G}, we have to work over the field $\RM$.
An eigenvalue of $\Delta_d$ can then be written as a vector $\vec{\lambda}_\gamma \in V$.
The action of the Cuntz-Krieger algebra $\Oo_{\tilde{A}}$ on such an eigenvalue, given in equation \eqref{tLap09.eq-ui(eval)}, becomes here
\begin{equation}
\label{tLap09.eq-ui(eval)-ring}
u_e \bigl( \vec{\lambda}_\gamma \bigr) \; =\; C \,\vec{\lambda}_\gamma  \; + \; \vec{\beta}_e\,,
\end{equation}
where $\vec{\beta}_e$ is the expression of $\beta_e$ as a vector in $V$.
The general expression of an eigenvalue given in proposition \ref{tLap09.prop-ui(eval)-explicit}, equation \eqref{tLap09.eq-CK-spectrum-eval}, takes here the form
\begin{equation}
\label{tLap09.eq-CK-spectrum-eval-ring}
\vec{\lambda}_\gamma \; = \; C^n \, \vec{\lambda}_{\varepsilon'} \; + \;
\sum_{j=1}^{n} \, C^{j-1} \, \vec{\beta}_{e_j} \,.
\end{equation}

We now assume that $C$ is {\em Pisot}, {\it i.e.} it has a single real eigenvalue $\mu_u >1$, and all its other eigenvalues $\mu$ satisfy $|\mu| < 1$.
Let us denote by $V_u$ the eigenspace of $\mu_u$, and by $V_s$ its orthogonal complement in $V$.
We can now characterize the distribution of the eigenvalues of $\Delta_d$.
\begin{theo}
\label{tLap09.thm-dist-eval}
There exist a constant $m>0$ such that for any eigenvalue $\lambda_\gamma$ of $\Delta_d$, one has
\begin{equation}
\label{tLap09.eq-dist-eval}
\text{\rm dist}_V \bigl( \vec{\lambda}_\gamma , V_u \bigr) \le \frac{m}{1 - \| C_s\|}\,,
\end{equation}
where $C_s$ is the projection of $C$ to $V_s$, and $\text{\rm dist}_V$ the induced Euclidean distance on $V$.
\end{theo}
\begin{proof}
The distance \(\text{\rm dist}_V \bigl( \vec{\lambda}_\gamma , V_u \bigr)\) equals the norm of the projection of $\vec{\lambda}_\gamma$ to $V_s$:
\begin{eqnarray*}
\text{\rm dist}_V \bigl( \vec{\lambda}_\gamma , V_u \bigr)  = &
\bigl\| 
P^{-1} \bigl(
C_s^n P^{-1} \vec{\lambda}_{\varepsilon'} + \sum_{j=1}^n C_s^{j-1} P^{-1} \vec{\beta}_{e_j} 
\bigr)
\bigr\| \\
\le & \|P^{-1}\|^2 \
\max \{ \|\vec{\beta}_e\|, \| \vec{\lambda}_{\varepsilon}\| \, : \, \varepsilon\in\Ee_0, e \in \Ee\}
\ \sum_{j=0}^{n}\|C_s^j\| \\
\le & m \; \sum_{j=0}^{\infty} \; \|C_s\|^j \; = \; m / (1-\|C_s\|)\,,
\end{eqnarray*}
where $P$ is an invertible matrix diagonalizing $C$, \(m = \|P^{-1}\|^2\max \{ \|\vec{\beta}_e\|, \| \vec{\lambda}_{\varepsilon}\| \}\), and the series converges because $C_s$ is diagonal with eigenvalues of moduli strictly less than one.
\end{proof}

\begin{rem}
\label{tLap09.rem-dist-eval-NPisot}
{\em If $C$ is no longer Pisot but {\em strictly hyperbolic} ({\it i.e.} has no eigenvalue of modulus 1), then the above result still holds but with $V_u$ the {\em unstable space} (the span of the eigenvectors with eigenvalues $|\mu| >1$) and $V_s$ its orthogonal complement.
}
\end{rem}

\section{Examples}
\label{tLap09.sect-Ex}

We illustrate here the results of sections \ref{tLap09.sect-LBop} and \ref{tLap09.sect-CK} for the classic examples of the Thue--Morse, Fibonacci, Ammann--A2, and Penrose tilings.

\subsection{The Fibonacci diagram}
\label{tLap09.ssect-fibo}

The Bratteli diagram for the (uncollared) Fibonacci substitution 
\(\left\{ \begin{array}{ccl} a & \rightarrow & ab \\ b & \rightarrow & a \end{array} \right.\)
reads
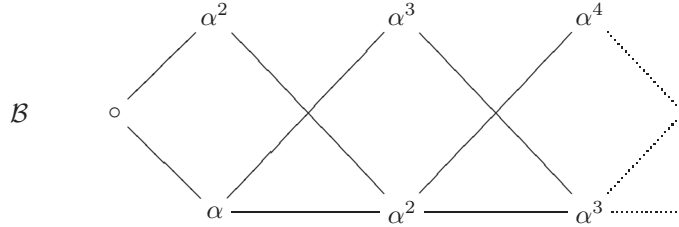
\[
\xymatrix{
&   &  \alpha^2 \ar@{-}[rrdd] & &  \alpha^3 \ar@{-}[rrdd] & &  \alpha^4 \ar@{.}[dr] & \\
\Bb &\circ \ar@{-}[ur] \ar@{-}[dr] & & & & & & \\
&   &  \alpha \ar@{-}[rr]  \ar@{-}[uurr] & &  \alpha^2\ar@{-}[rr] \ar@{-}[uurr] & &  \alpha^3 \ar@{.}[r] \ar@{.}[ur] & \\
}
\]
where the top vertices are of ``type $b$'' and the bottom ones of ``type $a$'', where \(\alpha= 1 /\phi=(\sqrt{5}-1)/2\) is the inverse of the golden mean, and the term $\alpha^n$ at a vertex is the measure of the cylinder of infinite paths through that vertex.

Note that this substitution does not force the border, so that $\partial \Bb$ is not the transversal of the Fibonacci tiling space.
For illustration purposes it is however worth carrying this example in details.
We treat the ``real'' Fibonacci tiling together with the Penrose tiling in section \ref{tLap09.ssect-PAF}.

Since the Bratteli diagram has only simple edges, as noted in Remark \ref{tLap09.rem-simple-edges}, the paths can be indexed by the vertices they go through.
The paths in $\Pi_1$ are thus written $a$ and $b$, and the paths in $\Pi_2$ are written $aa$, $ab$, and $ba$
(note that these are {\em not} orthonormal bases for the dot product given by the Dixmier trace , so that the Laplacians written below will not be symmetric).
The restrictions of the Laplace operator \eqref{tLap09.eq-LapExplicit} for $s=s_0 =1$ to $\Pi_1$ and $\Pi_2$ are given below together with their eigenelements.
\[
\Delta \vert_{\Pi_1} =
\left[
\begin{array}{cc}
1-\phi^2 & -1 + \phi^2 \\
\phi^2 & -\phi^2
\end{array}
\right]\,,
\ \text{\rm with eigenelements} \ \
(
0 \,, \, 
\left[
\begin{array}{c}
1 \\
1
\end{array}
\right] 
) \,,
\quad
(
1 + 2\phi^2 \,, \,
\left[
\begin{array}{c}
1 \\
1-\phi^2
\end{array}
\right] 
)\,.
\]
\begin{multline*}
 \Delta \vert_{\Pi_2} =
\left[
\begin{array}{ccc}
2 - 3\phi^2 & -1 +2 \phi^2 &  -1 + \phi^2 \\
-1 +3\phi^2 & 2 - 4\phi^2  & -1 +\phi^2\\
-1 + \phi^2 & 1 & -\phi^2
\end{array}
\right]\,,
\quad \quad \text{\rm with eigenelements} \\
\quad \quad
(
0, 
\left[
\begin{array}{c}
1 \\
1 \\
1
\end{array}
\right]
) \,, \
(
1 + 2\phi^2,
\left[
\begin{array}{c}
1 \\
1 \\
1-\phi^2
\end{array}
\right] 
) \,, \
(
3 + 6\phi^2,
\left[
\begin{array}{c}
1 \\
1-\phi^2\\
0
\end{array}
\right]
) \,.
\end{multline*}
Using the identities \(\chi_{b} = \chi_{ba}\), and \(\chi_a = \chi_{aa} + \chi_{ab}\), we see that the first two eigenvectors of $\Delta\vert_{\Pi_2}$ are exactly those of $\Delta\vert_{\Pi_1}$ expressed in $\Pi_2$.
Note that \( 1-\phi^2 = -\phi\) and that the above eigenvectors of $\Delta$ are \(\chi_{a} + \chi_b\), \(\chi_a - \phi \chi_b \), and \(\chi_{a a} - \phi \chi_{a b} \).
And all other eigenvectors are given by \(\chi_{\gamma a a} - \phi \chi_{\gamma a b} \) for $\gamma \in \Pi$.

Since the Bratteli diagram has only simple edges, as noted in Remark \ref{tLap09.rem-vertices}, we can take \( \tilde{\Bb} = \Bb\) and 
\( \tilde{A}= A = \left[ \begin{array}{cc} 1 & 1 \\ 1 & 0 \end{array} \right]\).
The action of the two Cuntz--Krieger operators $U_a$ and $U_b$ on the eigenvalues of $\Delta$ as in equation \eqref{tLap09.eq-actionCK-eval} is given in here by
\begin{equation}
\label{tLap09.eq-CKFibo}
u_a \bigl(\lambda_\gamma \bigr)= \phi^2 \lambda_\gamma + 1 -\phi^2\,, \quad \quad 
u_b \bigl( \lambda_\gamma \bigr) = \phi^2 \lambda_\gamma - 1 + \phi^2 \,, 
\end{equation}
if $\gamma$ is compatible with their action, and $u_a (\lambda_\gamma) = 0$ or  $u_b (\lambda_\gamma) = 0$ otherwise.

Over the ring \(\ZM \oplus \phi^2 \ZM\), the companion matrix of $A$ is \(\left[ \begin{array}{cc} 0 & 1 \\ 1 & 1 \end{array} \right]\), and the operators \eqref{tLap09.eq-CKFibo} become the affine maps
\[
\label{tLap09.eq-CKfibo}
u_a \bigl( \vec{\lambda}_\gamma \bigr) =  
\left[ \begin{array}{cc} 0 & 1 \\ 1 & 1 \end{array} \right] \vec{\lambda}_\gamma 
- \left[ \begin{array}{r} -1 \\ 1 \end{array} \right]\,, \quad \quad 
u_b \bigl( \vec{\lambda}_\gamma \bigr) =  
\left[ \begin{array}{cc} 0 & 1 \\ 1 & 1 \end{array} \right] \vec{\lambda}_\gamma 
+ \left[ \begin{array}{r} -1 \\ 1 \end{array} \right]\,, 
\]
when  $\gamma$ is compatible with the corresponding action.

Figure \ref{tLap09.fig-fibo} illustrates Theorem \ref{tLap09.thm-dist-eval} that characterizes the repartition of the eigenvalues of $-\Delta$ as point of integer coordinates that stay within a bounded strip to the Perron-Froebenius eigenline of $A$ (slope $\phi$) in \(\ZM \oplus \phi^2 \ZM\).
Note that the repartition of points in the strip is not ``homogeneous'', {\it i.e.} the number of points within a distance $r$ to the origin is not linear in $r$, but rather follows the Weyl asymptotics in $\sqrt{r}$ (Theorem \ref{tLap09.thm-Weyl}).

\subsection{The dyadic  Cantor set and the Thue--Morse tiling}
\label{tLap09.ssect-ThueMorse}

Those examples have enough symmetries to allow easy and direct calculations (without using the operators of the Cuntz--Krieger algebra).
The Bratteli diagram $\Bb$ of the dyadic Cantor set is the dyadic odometer, 
\[
\xymatrix{
\Bb  &  \circ   \ar@{=}[rr]^0_1 & & \frac{1}{2}  \ar@{=}[rr]^0_1 & &   \frac{1}{4}  \ar@{=}[rr]^0_1 & &  \frac{1}{8}  \ar@{:}[r]&
}
\]
and its associated diagram for its Cuntz-Krieger algebra $\tilde{\Bb}$ is the Bratteli diagram of the (uncollared) Thue-Morse substitution 
\(\left\{ \begin{array}{ccl} 0 & \rightarrow & 01 \\ 1 & \rightarrow & 10 \end{array} \right.\):
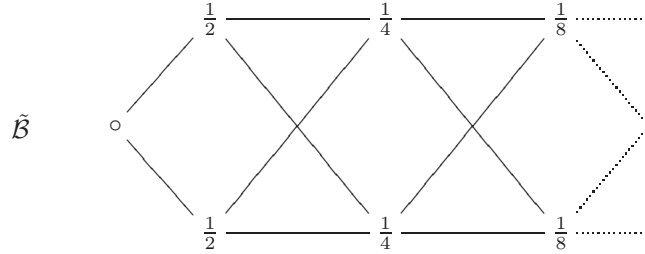
\[
\xymatrix{
&   &  \frac{1}{2} \ar@{-}[rrdd]  \ar@{-}[rr] & &   \frac{1}{4}  \ar@{-}[rrdd] \ar@{-}[rr] & &  \frac{1}{8}  \ar@{.}[dr] \ar@{.}[r]& \\
\tilde{\Bb} &\circ \ar@{-}[ur] \ar@{-}[dr] & & & & & & \\
&   &   \frac{1}{2} \ar@{-}[rr]  \ar@{-}[uurr] & &   \frac{1}{4} \ar@{-}[rr] \ar@{-}[uurr] & &   \frac{1}{8}  \ar@{.}[r] \ar@{.}[ur] & \\
}
\]
where the term $\frac{1}{2^n}$ at a vertex is the measure of the cylinder of infinite paths through that vertex.
The top vertices are of ``type $0$'', the bottom ones of ``type $1$''.

We label the paths in $\Bb$ by sequences of $0$'s and $1$'s labeling the edges they go through from the root: \( \gamma \in \Pi_n\) is written \( \gamma =(\eps_1, \cdots, \eps_n)\).
The Laplacian on $\Bb$ commutes with the following operators:
\begin{equation}
\label{tLap09.eq-taui}
\tau_i \chi_{(\eps_1, \cdots, \eps_n)} = 
\left\{
\begin{array}{ll}
\chi_{(\eps_1, \cdots, \eps_i + 1, \cdots \eps_n)} & \text{\rm if} \, 1 \le i \le n\,,\\
\chi_{(\eps_1, \cdots, \eps_n)} &  \text{\rm else,}
\end{array}
\right.
\end{equation}
where the addition is taken$\mod 2$.
The operators $\tau_i$ commute with each other and square up to the identity.
One can therefore choose an eigenbasis for $\Delta$ made of eigenelements of the $\tau_i$'s: that is {\em Haar functions} on the dyadic Cantor set $\partial \Bb$.
We recover this way the example treated in \cite{PB09} and we refer the reader there for the details.

For the Thue--Morse diagram $\tilde{\Bb}$, we can also index paths by sequences of $0$'s and $1$'s labeling the vertices they go through from the root.
The Laplacian is also commuting with the operators $\tilde{\tau}_i$ defined like the $\tau_i$ defined in equation \eqref{tLap09.eq-taui}.
A basis of eigenvectors of $\Delta$ for $s=s_0=1$ is given by the constant function \(\chi_{\partial \tilde{\Bb}}\) (with eigenvalue $0$), and the functions
\[
\varphi_{n, \gamma} = \chi_\gamma - \tilde{\tau}_n \chi_\gamma \,, \quad \gamma \in \tilde{\Pi}_n \,,
\]
for $n\in \NM$, with eigenvalues \(\lambda_n = -\frac{2}{3}\bigl( 7 \cdot 4^{n-1} - 1 \bigr)\) of degeneracy \(2^{n-1} = \text{\rm Card} \, \Pi_n\).
The eigenvalues satisfy the induction formula \(\lambda_{n+1} = 4 \lambda_n -2\).

The Weyl asymptotics of Theorem \ref{tLap09.thm-Weyl} reads here
\( \frac{1}{2} \sqrt{\frac{6}{7} \lambda + \frac{10}{7}} \le \Nn(\lambda) \le \sqrt{\frac{6}{7} \lambda + \frac{4}{7}}\).

\subsection{The Penrose tiling}
\label{tLap09.ssect-PAF}

The Fibonacci, Penrose, and Ammann--A2 \cite{GS89} tilings have formally the ``same'' substitution on prototiles modulo their symmetry groups, with Abelianization matrix
\[
A =
\left[ 
\begin{array}{cc}
2 & 1\\
1 & 1
\end{array}
\right] \,.
\]
The Penrose and Ammann--A2 substitutions force the border.
And for the Fibonacci tiling, one considers the conjugate substitution \(a\rightarrow baa, b \rightarrow ba\), which is primitive, recognizable, and forces the border as noted in Example \ref{tLap09.ex-fiboconj}.
Those three substitution tilings have the same Perron--Frobenius eigenvalue, namely $\pf = \phi^2$, where \(\phi = (1+\sqrt{5})/2\) is the golden mean. 
In conclusion, the transversals of those tiling spaces can be described by the set of infinite paths in the same Bratteli diagram $\Bb$ illustrated below:
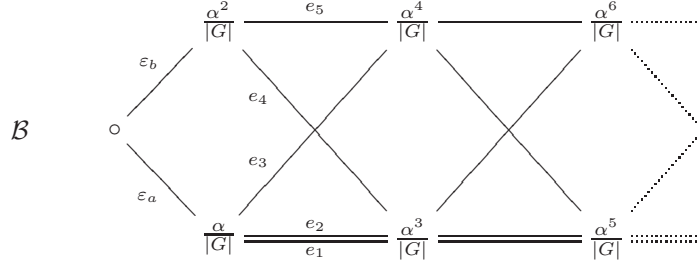
\[
\xymatrix{
&   &  \frac{\alpha^2}{|G|} \ar@{-}[rr]^{e_5} \ar@{-}[rrdd]_(.3){e_4} & &  \frac{\alpha^4}{|G|}\ar@{-}[rr] \ar@{-}[rrdd] & &  \frac{\alpha^6}{|G|} \ar@{.}[r] \ar@{.}[dr] & \\
\Bb &\circ \ar@{-}[ur]^{\varepsilon_b} \ar@{-}[dr]_{\varepsilon_a}& & & & & & \\
&   &  \frac{\alpha}{|G|} \ar@{=}[rr]_{e_1}^{e_2}  \ar@{-}[uurr]^(.3){e_3} & &  \frac{\alpha^3}{|G|}\ar@{=}[rr] \ar@{-}[uurr] & &  \frac{\alpha^5}{|G|} \ar@{:}[r] \ar@{.}[ur] & \\
}
\]
where \(\alpha= 1/\phi=(\sqrt{5}-1)/2\) is the inverse of the golden mean,
and the term $\alpha^n/|G|$ at a vertex is the measure of the set of infinite paths through that vertex.
And where $G$ is the symmetry group of the tiling introduced in section \ref{tLap09.rem-group}, and $|G|$ the cardinality of $G$.
That is $G=\{1\}$ is the trivial group (so $|G|=1$) for the Fibonacci tiling, $G=C_2 \times C_2$ for the Ammann--A2 tiling (symmetries of ``vertical and horizontal'' reflections, $|G|=4$), and $G=D_{10}$ for the Penrose tiling (10-fold rotational symmetries, and reflections, $|G|=20$).

Let us denote by $U_i, i= 1,\cdots 5,$ the generators of the Cuntz--Krieger algebra \eqref{tLap09.eq-CKalg} associated with
the Abelianization matrix of $\tilde{\Bb}$.
The induced action on the eigenvalues of $\Delta = \Delta_{s_0}$ as in equation \eqref{tLap09.eq-actionCK-eval} reads here for Penrose and Ammann--A2 (for Fibonacci, $\pf^{2/d} = \phi^4$ has to replace $\phi^2$ in the following equations):

\[
\begin{split}
u_1 \bigl( \lambda_\gamma \bigr) & =  \phi^2 \lambda_\gamma + 
(1-\phi^2) \frac{\mu[\varepsilon_a]-\mu(\circ)}{G(\circ)} + 
\frac{\mu[\varepsilon_ae_1] - \mu[a]}{G(\varepsilon_a)} \\ 
u_2 \bigl( \lambda_\gamma \bigr) & =  \phi^2 \lambda_\gamma + 
(1-\phi^2) \frac{\mu[\varepsilon_a]-\mu(\circ)}{G(\circ)} + 
\frac{\mu[\varepsilon_ae_2] - \mu[a]}{G(\varepsilon_a)} \\
u_3 \bigl( \lambda_\gamma \bigr) & =  \phi^2 \lambda_\gamma + 
-\phi^2 \frac{\mu[\varepsilon_a]-\mu(\circ)}{G(\circ)} + 
\frac{\mu[\varepsilon_b]-\mu(\circ)}{G(\circ)} + 
\frac{\mu[\varepsilon_be_3] - \mu[a]}{G(\varepsilon_b)} \\
u_4 \bigl( \lambda_\gamma \bigr) & =  \phi^2 \lambda_\gamma + 
-\phi^2 \frac{\mu[\varepsilon_b]-\mu(\circ)}{G(\circ)} + 
\frac{\mu[\varepsilon_a]-\mu(\circ)}{G(\circ)} +
\frac{\mu[\varepsilon_ae_4] - \mu[a]}{G(\varepsilon_a)} \\
u_5 \bigl( \lambda_\gamma \bigr) & =  \phi^2 \lambda_\gamma + 
(1-\phi^2) \frac{\mu[\varepsilon_b]-\mu(\circ)}{G(\circ)} + 
\frac{\mu[\varepsilon_be_5] - \mu[\varepsilon_b]}{G(\varepsilon_b)} \\
\end{split}
\]
if $\lambda_\gamma$ is compatible with the operators.
Here $G= G_{s_0}$ as in equation \eqref{tLap09.eq-G}, so we have
\(G(\circ) = |G| (|G|-1) \bigl( (\frac{\alpha}{|G|})^2 + (\frac{\alpha^2}{|G|})^2 \bigr) + |G|^2 \frac{\alpha}{|G|} \frac{\alpha^2}{|G|}\),
\( G(\eps_a) = 2 \bigl( (\frac{\alpha^2}{|G|})^2 + (\frac{\alpha^3}{|G|})^2 \bigr) + 4 \frac{\alpha^3}{|G|} \frac{\alpha^4}{|G|}\)
and \( G(\eps_b) = 2 \frac{\alpha^3}{|G|}\frac{\alpha^4}{|G|}\).
The eigenelements of $\Delta \vert_{\Pi_2}$ are $0$ for \(\chi_{\partial \Bb}\), $\lambda_0$ for \( \chi_{\eps_a} - \phi \chi_{\eps_b}\), $\lambda_{\eps_a}$ for \( \chi_{\eps_a e} - \phi \chi_{\eps_a f}\,, e, f \in \ext(\eps_a)\), and  $\lambda_{\eps_b}$ for \( \chi_{\eps_b e} - \phi \chi_{\eps_b f}\,, e, f \in \ext(\eps_b)\), where
\[
\lambda_0 = \frac{-2|G| \bigl( |G| + 1 - 4 \phi^2 \bigr)}{|G|^2 - 10|G|+5} \,, \quad
\lambda_{\eps_a}= \frac{\mu[\eps_a]-\mu[\circ]}{G(\circ)} - \frac{\mu[\eps_a]}{G(\eps_a)}\,, \quad
\lambda_{\eps_b}= \frac{\mu[\eps_b]-\mu[\circ]}{G(\circ)} - \frac{\mu[\eps_b]}{G(\eps_b)}\,.
\]

\bibliography{bibtLap09_12}
\bibliographystyle{plain}

\end{document}

%% file: repart5.pstex_t
\begin{picture}(0,0)%
\includegraphics{repart5.pstex}%
\end{picture}%
\setlength{\unitlength}{4144sp}%
\begingroup\makeatletter\ifx\SetFigFont\undefined%
\gdef\SetFigFont#1#2#3#4#5{%
  \reset@font\fontsize{#1}{#2pt}%
  \fontfamily{#3}\fontseries{#4}\fontshape{#5}%
  \selectfont}%
\fi\endgroup%
\begin{picture}(5615,3765)(-818,-4150)
\put(566,-3117){\makebox(0,0)[lb]{\smash{{\SetFigFont{10}{12.0}{\rmdefault}{\mddefault}{\updefault}{\color[rgb]{0,0,0}$u_a$}%
}}}}
\put(1757,-3117){\makebox(0,0)[lb]{\smash{{\SetFigFont{10}{12.0}{\rmdefault}{\mddefault}{\updefault}{\color[rgb]{0,0,0}$u_b$}%
}}}}
\put(1530,-2352){\makebox(0,0)[lb]{\smash{{\SetFigFont{10}{12.0}{\rmdefault}{\mddefault}{\updefault}{\color[rgb]{0,0,0}$u_a$}%
}}}}
\put(1644,-1501){\makebox(0,0)[lb]{\smash{{\SetFigFont{10}{12.0}{\rmdefault}{\mddefault}{\updefault}{\color[rgb]{0,0,0}$u_a$}%
}}}}
\put(2183,-1076){\makebox(0,0)[lb]{\smash{{\SetFigFont{10}{12.0}{\rmdefault}{\mddefault}{\updefault}{\color[rgb]{0,0,0}$u_a$}%
}}}}
\put(2977,-1813){\makebox(0,0)[lb]{\smash{{\SetFigFont{10}{12.0}{\rmdefault}{\mddefault}{\updefault}{\color[rgb]{0,0,0}$u_b$}%
}}}}
\put(3132,-720){\makebox(0,0)[lb]{\smash{{\SetFigFont{10}{12.0}{\rmdefault}{\mddefault}{\updefault}{\color[rgb]{0,0,0}$y=\phi x$}%
}}}}
\put(957,-517){\makebox(0,0)[lb]{\smash{{\SetFigFont{10}{12.0}{\rmdefault}{\mddefault}{\updefault}{\color[rgb]{0,0,0}$\phi^2 \mathbb{Z}$}%
}}}}
\put(3679,-3990){\makebox(0,0)[lb]{\smash{{\SetFigFont{10}{12.0}{\rmdefault}{\mddefault}{\updefault}{\color[rgb]{0,0,0}$\mathbb{Z}$}%
}}}}
\end{picture}%